\documentclass{article}
\usepackage{graphics}
\usepackage{amsfonts,amsmath}
\pdfoutput=1
\usepackage{amssymb}
\usepackage{subfigure}
\usepackage{comment}
\usepackage{cite}
\usepackage{float}
\usepackage{hyperref}
\usepackage{wrapfig}
\usepackage[dvipsnames]{xcolor}
\usepackage{graphicx} 
\usepackage{caption}
\usepackage{booktabs} 
\usepackage[T1]{fontenc} 
\usepackage{lmodern} 
\usepackage[utf8]{inputenc}

\usepackage{fancyhdr}
\usepackage[cm]{fullpage}
\usepackage{geometry}[left=0.5cm, right=0.5cm, top=0.95cm, bottom=0.95cm]

\newcommand\blfootnote[1]{%
  \begingroup
  \renewcommand\thefootnote{}\footnote{#1}%
  \addtocounter{footnote}{-1}%
  \endgroup
}

\title{Bifurcations of sleep patterns due to homeostatic and circadian variation in a sleep-wake flip-flop model}
\author{Christina Athanasouli \thanks{
Department of Mathematics, University of Michigan, Ann Arbor, MI (\hyperref[chrath@umich.edu]{chrath@umich.edu})} \and Sofia H. Piltz \thanks{Department of Mathematics, University of Michigan, Ann Arbor, MI (\hyperref[sofia.h.piltz@gmail.com]{sofia.h.piltz@gmail.com})} \and Cecilia Diniz Behn \thanks{Department of Applied Mathematics and Statistics, Colorado School of Mines, Golden, CO; Department of
Pediatrics, University of Colorado Anschutz Medical Campus, Aurora, CO (\hyperref[cdinizbe@mines.edu]{cdinizbe@mines.edu})} \and Victoria Booth \thanks{Departments of Mathematics and Anesthesiology, University of Michigan, Ann Arbor, MI (\hyperref[vbooth@umich.edu]{vbooth@umich.edu})}
}

\begin{document}

\maketitle

\begin{abstract}

Differential equation-based physiological models of sleep-wake networks describe sleep-wake regulation by simulating the activity of wake- and sleep-promoting neuronal populations and the modulation of these populations by homeostatic and circadian ($\sim24$ h) drives. Here, we consider a sleep-wake flip-flop network model consisting of mutually inhibitory interactions between wake- and sleep-promoting neuronal populations. Motivated by changes in sleep behavior during early childhood as babies transition from napping to non-napping behavior, we vary homeostatic and circadian modulation and analyze effects on resulting sleep-wake patterns. To identify the types and sequences of bifurcations leading to changes in stable sleep-wake patterns in this piecewise-smooth model, we employ multiple mathematical methods, including fast-slow decomposition and numerical computation of circle maps. We find that the average daily number of sleeps exhibits a period adding sequence as the homeostatic time constants are reduced, and that the temporal circadian profile influences the number of observed solutions in the sequence. These solutions emerge through sequences of saddle-node and border collision bifurcations, where the particular sequence depends on parameter values. When the temporal circadian profile is steep, as is observed with long day lengths, some sleep patterns are lost and bistability of other patterns may occur. We analyze a limiting case of the temporal circadian waveform, a circadian hard switch model, to understand this loss of solutions. Generally, our holistic analysis approach provides an alternative analysis method for model systems that defy conventional numerical bifurcation analysis techniques.

\end{abstract}
\blfootnote{\textbf{Funding:} This work was funded by the NSF under grants DMS 1853506 (C.A., S.H.P. and V.B.) and DMS 1853511 (C.D.B.).

\hspace{0.22cm}Preprint submitted to SIAM Journal on Applied Dynamical Systems (SIADS) in September 2021.}


\section{Introduction}
\label{Introduction}

Sleep timing is governed by interactions between circadian ($\sim 24$ h) and homeostatic sleep drives and the action of these drives on the networks of brainstem and hypothalamic neurons that promote states of wake and sleep \cite{Saper_2001, Saper2005}. Early models of sleep-wake behavior described this system with coupled oscillators representing interactions between sleep cycles driven by the homeostatic sleep drive and a periodic circadian rhythm \cite{Kronauer1982,Winfree1982,Daan1984,MooreEde1984,Winfree1983,Strogatzsleep}. The classic Two Process model quantified this relationship in the form of a threshold system consisting of an exponentially increasing and decreasing homeostatic sleep drive that switches direction at thresholds modulated by the periodic circadian rhythm \cite{tpm,Daan1984,Borbelyetal2016}. 
 
Circle maps have long been employed to study various biological systems consisting of coupled oscillators in which one oscillator drives another \cite{ArnoldArrhythmias, Glass1991, skeldon_dirks2019}. In addition, there exists a vast literature  on the analysis of functions of the circle to itself describing circle maps (e.g., \cite{KatokHasselblatt, GlassMackey,Keener1980,Boyland1986,Glendinning1995}). As reduced models for coupled oscillator systems, these results provide a powerful framework for understanding the dynamics of relative frequencies of coupled oscillators, including identifying types of phase-locked or entrained solutions \cite{HoppKeen,KeenHoppRinz,BaesensMacKay}, bifurcations between these solutions \cite{Nakao1998,skeldonphilips} and chaotic dynamics \cite{Glassetal,Guevaraetal,Ostbornetal}. Circle maps can be explicitly formulated for some model systems, such as threshold systems  \cite{Glass1991,Nakao1998,BRESSLOFF1990187,skeldon_dirks2019,Glendinning1995} or integrate-and-fire models \cite{ArnoldArrhythmias,skeldongaps}.
For example, circle maps of the Two Process model for the timing of sleep onset relative to the circadian rhythm can be explicitly computed, and recent work has analyzed these circle maps under several conditions \cite{skeldonphilips,skeldongaps}. These circle maps can be monotonic or non-monotonic, and may exhibit a vertical discontinuity or gap which introduces additional types of bifurcations between phase-locked solutions compared to continuous circle maps \cite{granados,AvrutinBook,DiBernardoBook}. 

However, using circle maps to understand solutions and their bifurcations in high- dimensional, differential equations-based, coupled oscillator models of biological processes is generally difficult since explicit computation of an underlying circle map is not straightforward. In this study, we numerically compute circle maps to identify the types of bifurcations in a model of a sleep-wake network under homeostatic and circadian
variation.

Physiologically-based models of sleep-wake networks are based on the interactions of neuronal populations that promote wake and sleep states, with the suprachiasmatic nucleus (SCN) that generates the circadian rhythm \cite{phillipsrobinson,DinizBehn2007,DBandBooth,GDBB,Rempe,KumarBose,BandDBMathBiosci}. The simplest of such ordinary differential equation (ODE) -based models consists of mutually inhibitory interactions between wake- and sleep-promoting populations, i.e.,  a sleep-wake flip-flop, with transitions dictated by homeostatic sleep and circadian rhythm drives \cite{phillipsrobinson}. One such model can be formally reduced to the same form as the Two Process model \cite{skeldonphilips} and a numerical study suggests that similar types of phase-locked solutions are obtained as the time constants of the homeostatic sleep drive and amplitude of circadian drive are varied \cite{PhillipsMammalianSleep}. However, to our knowledge, the types of bifurcations governing the gain and loss of stability of the phase-locked solutions in the ODE model have not been reported.

Here, we draw on multiple mathematical methods, including computation of circle maps, to analyze the bifurcations of phase-locked solutions in an ODE-based sleep-wake flip-flop model. In this model, interactions between neuronal populations are modulated by both homeostatic sleep and circadian drives. Motivated by changes in these drives that have been documented in the experimental literature, we consider the effects of variation in both the homeostatic sleep and the circadian drives on the types of bifurcations leading to changes in the number of daily sleep episodes. To our knowledge, this has not been thoroughly analyzed in a physiologically-based sleep-wake model previously. 

Specifically, we consider varying the time constants of the homeostatic sleep drive and the temporal profile of the circadian drive. Several experimental studies have shown that characteristics of sleep homeostasis depend on individual traits such as age and sex \cite{DijkSexDiffs, rusterholz,JenniAdolescents} and may vary with development \cite{rusterholzmain}. For example, the transition from polyphasic (multiple sleeps per day) to monophasic (one sleep per day) sleep behavior that occurs in infancy/early childhood is thought to result from differences in the time constants of the homeostatic sleep drive dictating the accumulation and dissipation of sleep pressure \cite{JENNI2007321,SALZARULO1992107, Jenni2004}. The temporal profile of the circadian drive, reflecting the firing rate of neurons in the SCN, can also vary with age \cite{CZEISLER1992933}, as well as in response to seasonal changes in day length. In particular, the mean duration of high SCN firing activity was shorter in animals entrained to a short photoperiod (i.e., daily illumination) (Light:Dark 8 h:16 h) and longer in those entrained to a long photoperiod (Light:Dark 14 h:10 h)\cite{Mrugala}. Similar results were reported in  \cite{Meijer, VANDERLEEST} who found longer intervals of peak firing activity of the SCN during long days and shorter peak firing intervals during short days. 

Previous analysis of mathematical models of sleep-wake regulation has considered the effects of changing homeostatic time constants in the transition from polyphasic to monophasic sleep \cite{PhillipsMammalianSleep,skeldonphilips,skeldongaps}. However, these studies have not taken into account the effects of the steepness of the circadian waveform and how this interacts with changing homeostatic time constants. Therefore, we extend these previous findings by additionally considering how the temporal profile of the circadian drive affects the bifurcation sequences of entrained sleep-wake patterns that are observed in the transition from polyphasic to monophasic sleep behavior and driven by changing homeostatic time constants.

To identify types of bifurcations and understand how they arise in the sleep-wake flip-flop (SWFF) model, we employ multiple methods to analyze model dynamics and numerically compile a two-parameter bifurcation diagram. First, we show that the model is a piecewise-smooth dynamical system of Filippov type with one switching boundary. Second, we take advantage of the separation of time scales between neuronal population activity (fast) and the homeostatic sleep and circadian drives (slow) to reveal stable and unstable manifolds that dictate the trajectory flows. Third, based on these fast-slow manifolds, we define and numerically compute circle maps for the timing of sleep onsets relative to the phase of the circadian rhythm. 

Similarly to the maps explicitly derived from the Two Process model, SWFF circle maps are discontinuous with an infinite slope on one side of the vertical discontinuity. Tracking fixed points on the circle maps as the homeostatic and circadian drives are varied allows us to characterize the types and sequences of bifurcations when different phase-locked solutions lose and gain stability. Finally, to explain changes observed with increasing steepness of the circadian drive waveform, we consider the limit of the circadian drive as a step function which we call the circadian hard switch model (CHS). This limit introduces a second switching boundary to the piecewise smooth model, and we analyze its bifurcation sequence as the homeostatic sleep drive is varied to verify the trends observed in the original SWFF model. Thus, with these multiple techniques and numerical simulations, we obtain a complete understanding of the dynamics of the SWFF model, and detect and classify the types of bifurcations that occur as two key parameters are varied. 


The paper is organized as follows: in Section \ref{Model_Tools} we introduce the sleep-wake flip-flop model and the different mathematical approaches that we employ to describe and analyze model dynamics. 
In Section \ref{Bifurcations} we analyze the bifurcations of entrained sleep-wake patterns under variation of the time constants of the homeostatic sleep drive and the duration of peak activity of the circadian drive waveform. In Section \ref{HSmodel} we formulate the circadian hard switch model and describe the bifurcations in this case. In Section \ref{Discussion}, we provide a brief discussion of our results.

\section{Sleep-Wake Flip-Flop (SWFF) model}
\label{Model_Tools}

\subsection{Model equations}

The canonical sleep-wake flip-flop (SWFF) model includes two neuronal populations that govern the transitions between the states of wake and sleep: a wake-promoting ($W$) and a sleep-promoting ($S$) population are coupled by mutual inhibition, and their interaction is modulated by homeostatic sleep and circadian drives. 
In our SWFF model, the circadian input is mediated by a third neuronal population representing the suprachiasmatic nucleus ($SCN$), a group of cells in the hypothalamus that acts as the circadian pacemaker  
and displays a 24-hour variation in neural firing. 
For humans under typical conditions, the circadian rhythm and the sleep-wake cycle are entrained with lower SCN firing rates during sleep and higher SCN firing rates during wake.

\begin{figure}[hbt!]
	\centering
	\includegraphics[width=15cm,height=5.1cm]{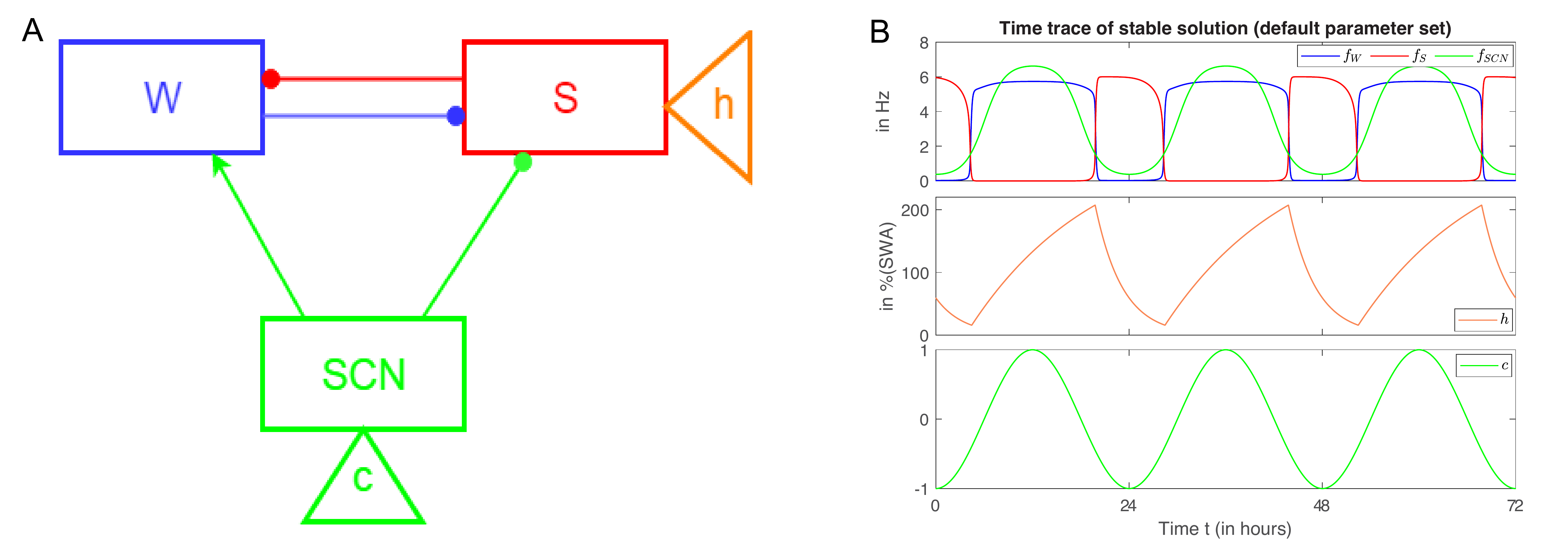}
	\caption{A SWFF model for sleep-wake regulation. A: Schematic of the model network summarizing interactions among the wake-promoting ($W$), sleep-promoting ($S$) and suprachiasmatic nucleus ($SCN$) neuronal populations with circles denoting inhibitory and arrows denoting excitatory synaptic connections. The homeostatic sleep drive ($h$) modulates activity of the sleep-promoting population and the circadian drive ($c$) modulates activity of both the sleep- and wake-promoting populations through SCN.  B: Time traces of the stable solution of the model for the default parameter set that resembles stereotypical adult human sleep. The firing rates for wake- ($f_W$, blue), sleep- ($f_S$, red) promoting and SCN ($f_{SCN}$, green) populations are shown in the top panel. The  middle and bottom panels include the profiles of the homeostatic sleep drive ($h$) and the circadian drive ($c$), respectively.}
	\label{fig:timetraces_default1}
\end{figure}

We use a firing rate formalism to model the neuronal population activity. Instead of tracking the spiking of single neurons, standard firing rate models describe the averaged behavior of spike rates of neuronal populations ($f_{W}$, $f_{S}$, $f_{SCN}$) \cite{wilsoncowan,Deco2008}. In particular, the mean postsynaptic firing rates are driven by the weighted mean firing rates of the presynaptic populations. 

\subsubsection*{Neuronal populations} The equations for the neuronal populations are as follows:

\begin{equation} \label{FF_fW}
\dfrac{df_{W}}{dt}=\dfrac{(W_{\infty}(g_{scnw}f_{SCN}-g_{sw} f_{S})-f_{W})}{\tau_{W}} \,,
\end{equation}
\begin{equation} \label{FF_fS}
\dfrac{df_{S}}{dt}=\dfrac{(S_{\infty}(-g_{ws} f_{W}-g_{scns} f_{SCN})-f_{S})}{\tau_{S}}\,,
\end{equation}
\begin{equation} \label{FF_fSCN}
\dfrac{df_{SCN}}{dt}=\dfrac{(SCN_{\infty}(c(t))-f_{SCN})}{\tau_{SCN}}\,.
\end{equation}

The postsynaptic firing rates, $f_{X}(t)$ (in Hz), saturate to their steady state firing rate response functions $X_{\infty}(\cdot)$ with time constants $\tau_{X}$ for $X= W, S, SCN$. The steady state firing rate functions, $X_{\infty}(\cdot)$, have a sigmoidal profile that has been utilized in many firing rate models \cite{wilsoncowan,onemap,phillipsrobinson, Deco2008}:


\begin{equation} \label{FF_Winf}
W_{\infty}(x)=W_{max}\cdot0.5\cdot\bigg(1+\tanh\Big(\dfrac{x-\beta_W}{\alpha_W}\Big)\bigg)\,,
\end{equation}
\begin{equation}\label{FF_Sinf}
S_{\infty}(x)=S_{max}\cdot0.5\cdot\bigg(1+\tanh\Big(\dfrac{x-\beta_S(h)}{\alpha_S}\Big)\bigg)\,,
\end{equation}
\begin{equation}\label{FF_fSCN}
SCN_{\infty}(x)=SCN_{max}\cdot0.5\cdot\left(1+\dfrac{\tanh\Big(\frac{1}{0.7}\Big)}{\tanh\Big(\frac{1}{\alpha_{SCN}}\Big)}\tanh\Big(\dfrac{x-\beta_{SCN}}{\alpha_{SCN}}\Big)\right)\,.
\end{equation}

\subsubsection*{Homeostatic sleep drive}
The homeostatic sleep drive ($h$) regulates sleep propensity and is based on experimentally observed variation in the power of slow wave (0.75 - 4.5 Hz) fluctuations in electroencephalogram (EEG) recordings during sleep \cite{rusterholzmain,Daan1984, borbely2000principles, tpm}. The levels of the homeostatic sleep drive increase exponentially with the time constant $\tau_{hw}$ while in wake 
and decrease exponentially with the time constant $\tau_{hs}$ during sleep according to
\begin{equation}\label{FF_h}
\dfrac{dh}{dt}=\dfrac{\mathcal{H}(f_{W}-\theta_{W})\cdot(h_{max}-h)}{\tau_{hw}}+\dfrac{\mathcal{H}(\theta_{W}-f_{W})\cdot(h_{min}-h)}{\tau_{hs}}\,,
\end{equation}
where $\mathcal{H}$ represents a Heaviside function and $h$ is in units of percent slow wave activity (SWA) power.  The time constants $\tau_{hw}$ and $\tau_{hs}$ are set to experimentally estimated values for typical adult human sleep behavior \cite{rusterholzmain}.
The sleep drive $h$ modulates the activity of the sleep-promoting population through the $h$-dependent activation threshold $\beta_{S}(h)$ as follows:
\begin{equation}\label{FF_betaS}
\beta_S(h)=k_2\cdot h+k_1\,.
\end{equation}
In this way as $h$ increases during wake, the sleep promoting population will activate to inhibit the wake population and cause the transition to sleep. Conversely, as $h$ decreases during sleep, the sleep population will inactivate and allow the wake population to activate. We define sleep onset to occur when $f_{W}$ decreases through $\theta_{W}$ (and $h$ starts to decrease) and wake onset to occur when $f_{W}$ increases through $\theta_{W}$ (and $h$ starts to increase).

\subsubsection*{Circadian drive}
The input to the SCN population is the circadian drive $c(t)$ which induces a 24h periodic variation of $f_{SCN}$. The input $c(t)$ is modeled by a simple sinusoidal function and is assumed to be entrained to the 24-hour day. 

	\begin{equation}\label{FF_cde}
\dfrac{dc}{dt}=-\omega\sin{\theta}\,,
\end{equation}
\begin{equation}\label{FF_theta}
\dfrac{d\theta}{dt}=\omega,\hspace{0.15cm} \textrm{where}\hspace{0.15cm}\omega=\dfrac{2\pi}{24}\,,
\end{equation}
which for an initial condition $\big(c(0),\theta(0)\big)=\big(\cos(-\phi \dfrac{2\pi}{24}), \phi\big)$  gives the stable solution:
\begin{equation}\label{FF_csol}
c(t)=\cos\Big((t-\phi)\cdot\dfrac{2\pi}{24}\Big)\,.
\end{equation}

\subsubsection*{Model parameters}
We have chosen our default parameter set (see Table \ref{tab:DefaultParameterSetTable}) to generate typical human sleep behavior similar to previous work \cite{onemap}. In Figure \ref{fig:timetraces_default1}B the time traces of the stable solution of the model are displayed. The wake and sleep durations, dictated by the time intervals when $f_W$ is above or below the threshold value $\theta_W$, respectively, are approximately 15.33 and 8.67 hours. As is typical for entrained adult human sleep, wake onset occurs at the early rise of the circadian cycle, while sleep onset occurs as SCN activity approaches its minimum. 

\begin{table}[H]
	\centering
	\begin{tabular}{||c|c|c|c||}
		\hline
		$W_{max}=6$ Hz & $\tau_W=0.1$ hr & $\alpha_W=0.5$ & $\beta_W=-0.37$\\
		\hline
		$S_{max}=6$ Hz & $\tau_S=0.1$ hr & $\alpha_S=0.175$ & $ $\\
		\hline
		$SCN_{max}=7$ Hz & $\tau_{SCN}=0.05$ hr & $\alpha_{SCN}=0.7$ & $\beta_{SCN}=0$\\
		\hline
		$g_{sw}=0.3$ & $g_{scnw}=0.06$ &  $g_{ws}=0.28$ & $g_{scns}=0.0825$ \\
		\hline
		$h_{max}=323.88$ & $h_{min}=0$ & $\tau_{hw}=15.78$ hr & $\tau_{hs}=3.37$ hr\\
		\hline
		$ k_1=-0.1$ & $k_2=-0.006$ & $\theta_W=4$ Hz & \\
		\hline
	\end{tabular}

\caption{Parameter values for the SWFF model. For $X	= W, S, SCN$, $\alpha_{X}$ and $\beta_{X}$ are in units of effective synaptic input. Additionally, for $Y= W, S$, $g_{XY}$ (where $X\neq Y$) has units of (effective synaptic input / Hz). Units for $h_{max}$ and $h_{min}$ are percentage mean SWA. The parameters $k_1$ and $k_2$ are measured in effective synaptic input and  effective synaptic input/($\%$ mean SWA), respectively. The remaining units are included in the table.} 
\label{tab:DefaultParameterSetTable}
\end{table}

\subsection{Summary of the model dynamics}

In this section, we analyze the model equations (Eq. \ref{FF_fW}-\ref{FF_theta}) with the default parameter values (see Table \ref{tab:DefaultParameterSetTable}) and introduce the relevant techniques employed to understand the model dynamics.

\subsubsection*{Piecewise smooth dynamical system}
Switching in the homeostatic sleep drive from increasing during wake to decreasing during sleep introduces a discontinuity in the derivative of $h$. On either side of this discontinuity, model dynamics are smooth, but the presence of the discontinuity can influence model trajectories at the boundary of the smooth regions.   In our system, the switching boundary is $\Gamma=\{f_{W}=\theta_{W}\}$, where $\theta_{W}=4$ Hz, that separates the system into two smooth vector fields.

	To formally define the model as a piecewise smooth system, let $\mathbf{X}= [ f_{W},f_{S},f_{SCN},h,c,\theta ]$. Define $\Gamma^{+}=\{f_{W}>\theta_{W}\}$ and $\Gamma^{-}=\{f_{W}<\theta_{W}\}$ as the regions on either side of $\Gamma$ where $F_{1}(\mathbf{X})$, $F_{2}(\mathbf{X})$, respectively, are the corresponding vector fields dictating model dynamics. Then we can rewrite our model system as follows:  

	\[ \dfrac{d\mathbf{X}}{dt}=\begin{cases} 
	F_{1}(\mathbf{X}) & \mathbf{X}\in \Gamma^{+}\\
	F_{2}(\mathbf{X}) & \mathbf{X}\in \Gamma^{-}\\
	\end{cases}
	\]

On $\Gamma^{+}$, $\dfrac{dh}{dt}=\dfrac{h_{max}-h}{\tau_{hw}}$ and on $\Gamma^{-}$, $\dfrac{dh}{dt}=\dfrac{h_{min}-h}{\tau_{hs}}$, while the rest of the differential equations are defined as above. Since the vector fields are discontinuous across the switching boundary $\Gamma$,  we have a Filippov system (see \cite{DiBernardoBook}).\\
Generally, in Filippov systems as model trajectories approach the switching boundary, they may move along or ``slide'' on the switching boundary depending on the directions of the vector fields on either side of the boundary. To ensure that sliding along $\Gamma$ does not occur in our system, we check that $\Gamma$ is never simultaneously attracting (or repelling) for the flows in the vector fields on both sides \cite{DiBernardoBook}. To that end, let us consider $g(\mathbf{X})=f_W-\theta_W=0$ to define the boundary $\Gamma$. Then, $\nabla g= [1,0,0,0,0,0 ] $ and on $\Gamma$:

		$$\big(\nabla g(\mathbf{X})^{T}\cdot F_{1}(\mathbf{X})\big)\big(\nabla g(\mathbf{X})^{T}\cdot F_{2}(\mathbf{X})\big)=$$
\begin{align}\label{FF_CrossingCondition}			&\Bigg(\dfrac{W_{max}\cdot0.5\cdot\Big(1+\tanh\Big(\dfrac{g_{scnw}f_{SCN}-g_{sw} f_{S}-\beta_{W}}{\alpha_W}\Big)\Big)-f_{W}}{\tau_{W}}\Bigg)^2\geq0
	\end{align}
for $\mathbf{X} = [\theta_W, f_{S},f_{SCN},h,c,\theta ]$.

Condition \ref{FF_CrossingCondition} implies that the directions of the vector fields at the switching boundary $\Gamma$ are the same on either side. Thus, model trajectories cross $\Gamma$ when transitioning from one vector field to the other and no sliding along $\Gamma$ occurs. 
It is important to note that for all $\mathbf{X}$ on $\Gamma$ where $\nabla g(\mathbf{X})^{T}\cdot F_{1}(\mathbf{X})$ is equal to zero, then $\nabla g(\mathbf{X})^{T}\cdot F_{2}(\mathbf{X})$ is also zero. Therefore, their product will remain positive even if both quantities change sign. 

\subsubsection*{Fast and slow subsystems in the model}
In our model, the homeostatic sleep drive $h$, the circadian input $c$, and the circadian phase $\theta$ vary much more slowly compared to the firing rates $f_{W}$, $f_{S}$ and $f_{SCN}$. Hence, there is a separation of time scales which allows us to divide our system into a fast and slow subsystem, consisting of neuronal firing rates ($f_W$, $f_S$, and $f_{SCN}$) and the variables $h$, $c$, and $\theta$, respectively. Following a  similar analysis as in \cite{fastslow,onemap,skeldonphilips}, 
we define $\tau=\max{\{\tau_{W},\tau_{S},\tau_{SCN}\}} << \min{\{\tau_{hw},\tau_{hs},1/\omega\}}=\chi$.  These time scales introduce the parameter $\epsilon=\frac{\tau}{\chi}$, where $\epsilon$ has small magnitude.
 Moreover, $\min{\{\tau_{W},\tau_{S},\tau_{SCN}\}}=\mu\tau$ and $\max{\{\tau_{hw},\tau_{hs}, 1/\omega\}}=\lambda\chi$, where $\lambda=O(1)$ and $\mu=O(1)$.  Let us call $\Tilde{t}=\frac{t}{\tau}$ the time variable of the fast subsystem and $T=\frac{t}{\chi}$ the time variable of the slow subsystem, such that $\frac{\Tilde{t}}{T}=\frac{1}{\epsilon}$. Making the change of variables in Equations (\ref{FF_fW} -- \ref{FF_h}) and (\ref{FF_cde}) leads to 
 \begin{equation}\label{FF_fastslow1}
     \dfrac{d\mathbf{X_{fast}}}{d\Tilde{t}}=\mathbf{M}(\mathbf{X_{fast}},\mathbf{X_{slow}})
 \end{equation}
\begin{equation}\label{FF_fastslow2}
     \dfrac{d\mathbf{X_{slow}}}{d\Tilde{t}}=\epsilon\mathbf{N}(\mathbf{X_{fast}},\mathbf{X_{slow}})
 \end{equation} 
 where $\mathbf{X_{fast}}=[f_{W}, f_{S}, f_{SCN}]$, $\mathbf{X_{slow}}=[h, c, \theta]$, $\mathbf{M}$ and $\mathbf{N}$ are mapped to $\mathbb{R}^3$. 
 Coupling between $\mathbf{X_{fast}}$ and $\mathbf{X_{slow}}$
 occurs due to the dependence of the steady state functions $S_{\infty}(\cdot)$ on $h$ and $SCN_{\infty}(\cdot)$ on $c$. We note that the separation of time scales is valid in both vector fields on either side of the switching boundary $\Gamma$. 

\begin{figure}[hbt!]
  \hspace{-0.5cm}  
  \centering
\includegraphics[scale=0.47]{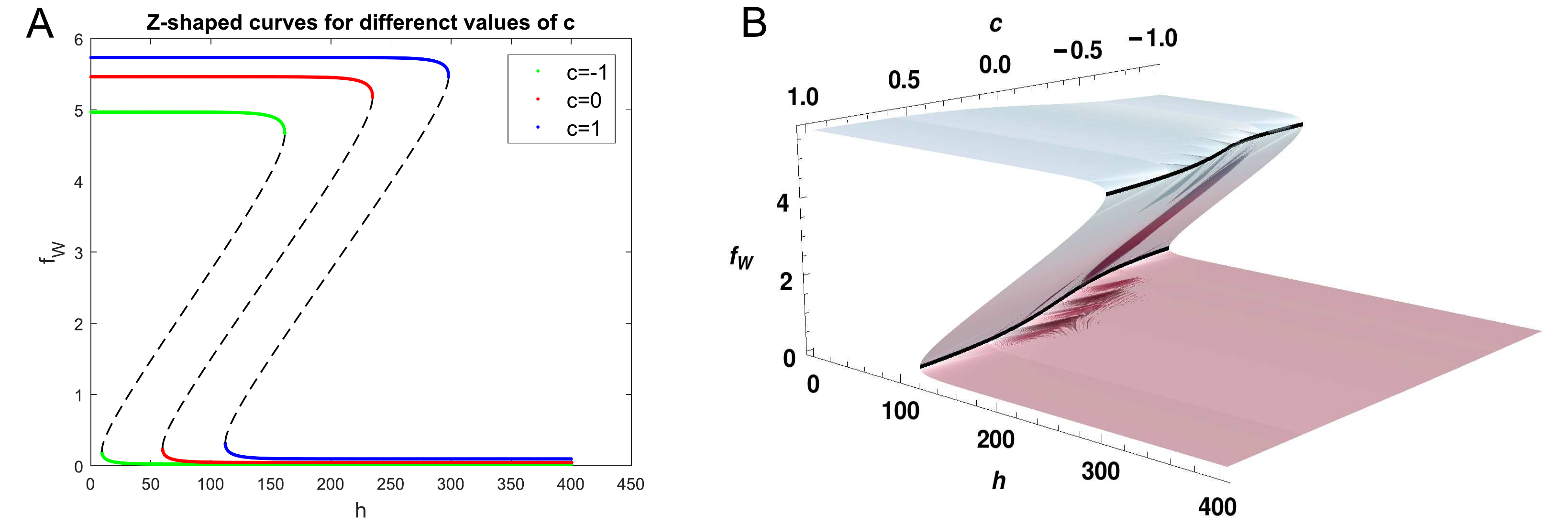}
\caption{Fast-slow decomposition of the SWFF model A: Bifurcation diagrams of the steady state solutions of the fast subsystem (Eq \ref{FF_fastslow1}) with respect to bifurcation parameter $h$ in the $h-f_W$ plane for fixed values of $c$ and $\theta$: $c=-1, \theta=12$ (green), $c=0, \theta=6$ (red) and $c=1, \theta=0$ (blue). The upper and lower branches (in green, red and blue) represent the stable steady wake and sleep states, respectively. The middle branch (dashed) represents an unstable state that separates the basins of attractions of the stable steady states. Notice that each curve has different upper and lower saddle-node points that eventually will define a curve in the $c-h-f_W$ space. B: Obtaining the $Z-$shaped curves for all values of $c\in [-1, 1]$ defines a surface that maintains the general $Z$-shape. The upper and lower saddle-node points of each $Z-$shaped curve define upper and lower saddle-node curves (curves in black). }
\label{fig:Zsurfacedefault}
\end{figure}

To analyze solution dynamics, we implement the fast-slow decomposition of our system by setting $\mathbf{X_{slow}} = \overline{\mathbf{X_{slow}}}$ to time-fixed parameters
and considering equilibrium solutions of the fast subsystem $\big(\dfrac{d\mathbf{X_{fast}}}{d\Tilde{t}} = 0\big)$. We represent the solutions of $\mathbf{M}(\mathbf{X_{fast}},\overline{\mathbf{X_{slow}}})=0$ 
in terms of the firing rate,
$f_{W}$, of the wake-promoting population. Specifically, setting $c$ and $\theta$ constant and computing solutions 
with respect to the bifurcation parameter $h$ yields 
a $Z$-shaped curve of steady states (Figure \ref{fig:Zsurfacedefault}A). The upper and lower branches of the $Z-$curve represent stable steady
states corresponding to the wake and sleep state, respectively. The middle
branch represents an unstable steady state that separates the basins of attraction of
the stable steady states. Finally,  the folds of each $Z-$curve are saddle-node bifurcation points
where the unstable steady state collides with one of the stable steady states.

For different (fixed) values of $c$ and $\theta$, the $Z-$shape of the curve is preserved, but the locations of the saddle-node bifurcation points with respect to $h$ change (Figure \ref{fig:Zsurfacedefault}A). Thus, as the circadian drive $c$ varies slowly, it affects the bifurcation structure of the fast subsystem. By definition, $c$ varies periodically between -1 and 1, and $\theta$ is such that $c=\cos(-\frac{\pi}{12}\theta)$. Hence, a $Z$-shaped surface is defined between these two extremes of the circadian cycle for the steady state solutions of the fast subsystem as a function of $h$ and $c$ (Figure \ref{fig:Zsurfacedefault}B).

When $h$, $c$ and $\theta$ vary slowly, model trajectories traverse the upper plane of the $Z-$surface during wake and the lower plane during sleep. Transitions between states occur when the trajectory reaches the curve of saddle-node points on either plane and evolves to the other plane. Note that the switching boundary $\Gamma$ lies between the upper (wake) and lower (sleep) planes of the $Z-$surface and trajectories  cross it during the transition. At $\Gamma$ crossing, $h$ changes direction leading to trajectory flows that follow a hysteresis loop around the $Z-$shaped surface. In this way, sleep (wake) onset is initiated when the model trajectory passes over the upper (lower) curve of saddle-nodes. 

Bifurcation diagrams were computed numerically using the software
AUTO XPPAUT \cite{Ermentrout}, and the fast-slow $Z$-shaped surfaces were created using Mathematica.

\vspace{0.5cm}

\subsubsection*{Sleep onset circle map}

 To analyze model dynamics and predict solution trajectories, we compute circle maps for the circadian phases of successive sleep onsets as in previous work \cite{onemap}. Specifically, define a Poincar{\'e} section for sleep onset as the firing rate of the wake-promoting population, $f_{W}$, decreasing through the switching boundary $\Gamma = \{ f_W = \theta_W \}$.  We define the circadian phase of the $n^{th}$ sleep onset, $\Phi_n$, as the time difference between the intersection of the model trajectory with the section (sleep onset)  and the preceding minimum of the $f_{SCN}$ variable  divided by the period of the circadian drive $c(t)$:
\begin{center}
     $\Phi_n$ =  $\dfrac{1}{24} (\textrm{time of sleep onset section crossing - time of preceding circadian minimum})$
\end{center}
We then define $\Pi: [0,1]\to [0,1]$ as the circle map with $\Phi_{n+1}=\Pi(\Phi_n)$. 

\begin{figure}[H]
\centering
	\includegraphics[scale=0.4]{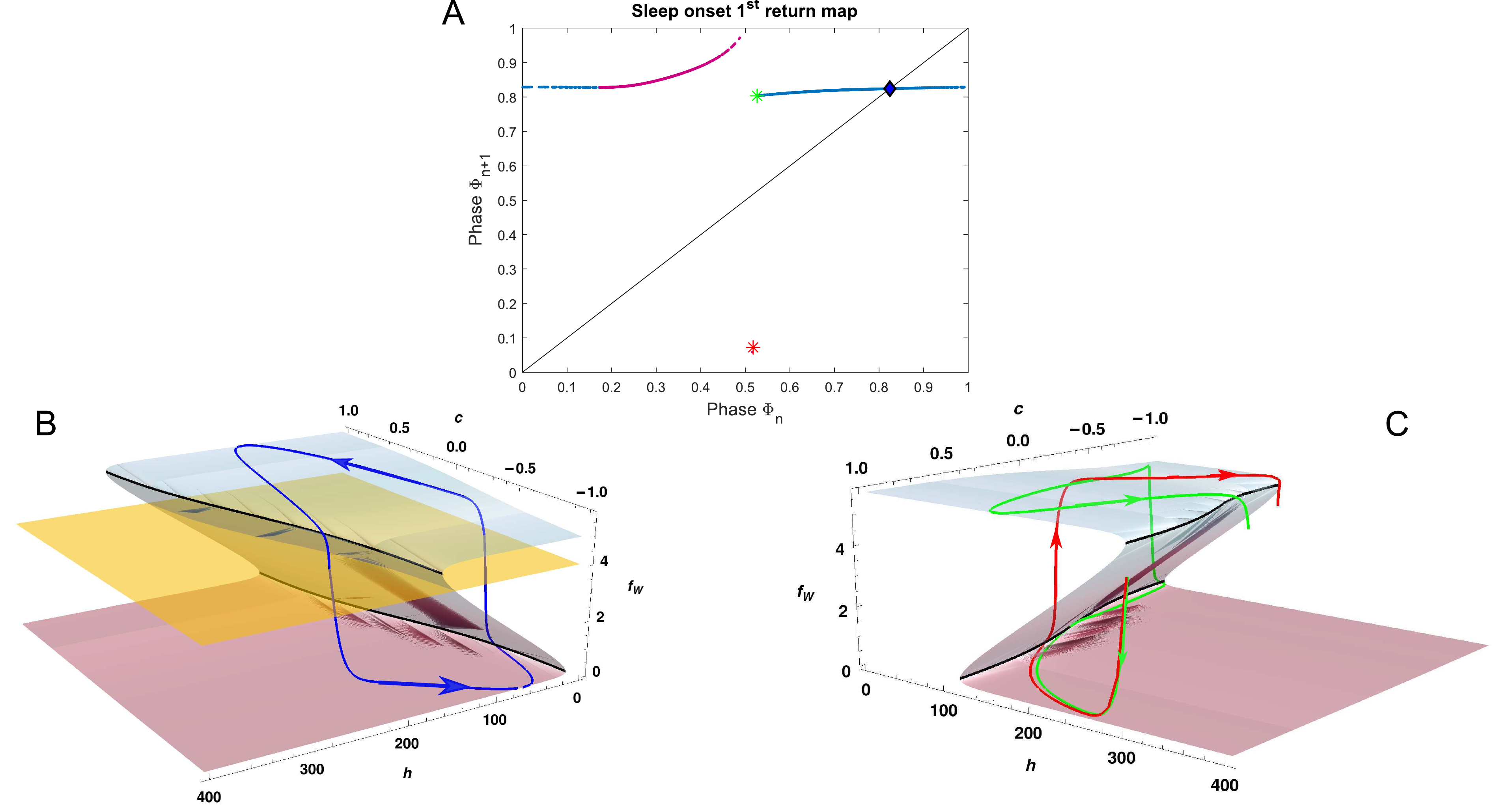}
\caption{Circle map and model trajectories relative to the fast-slow decomposition of the SWFF model. A: First return circle map for circadian phase of $n+1^{st}$ sleep onset, $\Phi_{n+1}$ as a function of circadian phase of $n^{th}$ sleep onset, $\Phi_n$. Purple points indicate circadian phases obtained by integrating the model from initial conditions ``forced" to lie on the unstable manifold. The blue diamond corresponds to the stable orbit shown in panel B (in blue). The green and red asterisks correspond to circadian phases associated with the trajectories in panel C (red and green, respectively).  B: The stable trajectory for the default parameter set (blue curve) plotted on the $Z-$surface computed from equilibrium solutions of the fast subsystem in Eq. (\ref{FF_fastslow1}). Sleep is initiated when the trajectory falls off the upper saddle-node curve. Sleep onset is defined as the time the trajectory crosses the switching boundary  $\Gamma = \{f_{W}=4\}$ (yellow plane) and $h$ starts decreasing. C: Trajectories with initial conditions on either side of the gap in the sleep onset map exhibit distinct behavior. The green trajectory becomes tangent to the lower saddle-node curve, resulting in a longer sleep episode, while the red one passes over the saddle-node curve and transitions to the wake state.}
\label{fig:defaultMAP_Trajectories}
\end{figure}

To compute the map $\Pi$, we simulate the model from initial conditions corresponding to sleep onset occurring at each circadian phase. 
Recall that the transition to sleep is initiated when the model trajectory passes over the curve of saddle-node points on the upper plane of the $Z$-shaped surface of steady state solutions of the fast subsystem. Thus, we select points on the 
upper saddle-node curve as a stable solution of the sleep-wake network that is near sleep
onset and use those values for initial conditions for the sleep-wake network variables in the map
construction. We compute these values on the upper saddle-node curve for all $c$ values over one circadian cycle by two-parameter numerical continuation, implemented in AUTO using XPPAUT \cite{Ermentrout}. 
By numerically integrating the model from these initial conditions, the circadian phases of sleep onsets are computed as the trajectories pass through the switching boundary $\Gamma$. 
From the majority of these initial conditions, model trajectories immediately transition down to the lower plane of the $Z-$shaped surface crossing $\Gamma$ on the way. However, 
 there is an interval of circadian phase values (where $\Phi_n$ is between approximately 0.2 and 0.4) for which the model trajectory does not immediately transition to sleep but instead continues to move along the upper plane until eventually transitioning at a later circadian phase. This produces a horizontal gap in the map.  This phenomenon occurs during the rising phase of $c$ that promotes the waking state at higher values of the homeostatic sleep drive $h$.  Dynamically speaking, the variables $c$ and $h$ vary on similar time scales at these moments, and the drive to sleep associated with increasing $h$ is balanced by a drive to wake associated with increasing $c$. To overcome this issue and fill in the horizontal gap in the map, for this interval of circadian phases, we substituted initial condition values that lie on the unstable manifold associated with the saddle of the upper saddle-node point that are closer to the switching boundary $\Gamma$ (purple points in Figure \ref{fig:defaultMAP_Trajectories}A).

Figure \ref{fig:defaultMAP_Trajectories}A illustrates the circle map; the circadian phase of the $n^{th}$ crossing of the section defined by $\Gamma$, $\Phi_{n}$, is on the $x-$axis, and the circadian phase of the  $n+1^{st}$ section crossing $\Phi_{n+1}$ is on the $y-$axis.
This first return sleep onset map is periodic in phase, consists of two branches, and has one stable fixed point at approximately (0.824,0.824), indicating that sleep onset of the stable solution occurs close to the trough of the circadian cycle. The stable solution trajectory is shown relative to the $Z-$shaped surface in Figure \ref{fig:defaultMAP_Trajectories}B. 
The map exhibits a discontinuity or gap around $\Phi_{n}=0.5$. The left branch of the discontinuity has an infinite slope which is a consequence of trajectories approaching a tangent intersection with the saddle-node curve of the $Z-$shaped surface (Figure \ref{fig:defaultMAP_Trajectories}C).  To see this, consider trajectories initiated on either side of the gap with sleep onsets very close to the peak of the circadian drive ($c=1$, red and green curves). The trajectory initiated on the infinite slope to the left of the gap (red curve) exhibits a short sleep episode, as it jumps up from the lower saddle-node curve and transitions to the wake plane resulting in the next sleep onset phase of about $\Phi_{n+1}=0.0722$. By contrast, the trajectory initiated on the right of the gap (green curve) becomes tangent to the lower saddle-node curve, resulting in a longer sleep episode. As the green trajectory evolves further, wake onset occurs close to the circadian minimum (that is $c=-1$), followed by a long wake episode resulting in the next sleep onset at a phase of about $\Phi_{n+1}=0.8033$.

    \section{Analysis of bifurcation sequences in the SWFF model}
    \label{Bifurcations}
    For our analysis of bifurcations in the SWFF model, we first identify the bifurcation sequences associated with the emergence and change of stable, phase-locked solutions as the time constants of the homeostatic sleep drive are varied. We then consider how the steepness of the circadian waveform affects the bifurcations of stable, phase-locked solutions and the associated bifurcation sequences as homeostatic time constants vary. We end this section by examining how the bifurcation sequences change in the regime of the fastest homeostatic time constants when sleep onset circle maps may be continuous.

    \subsection{Varying time constants of the homeostatic sleep drive}
  
\begin{figure}[H]
    	 \centering
\includegraphics[height=6.5cm, width = 15.4cm]{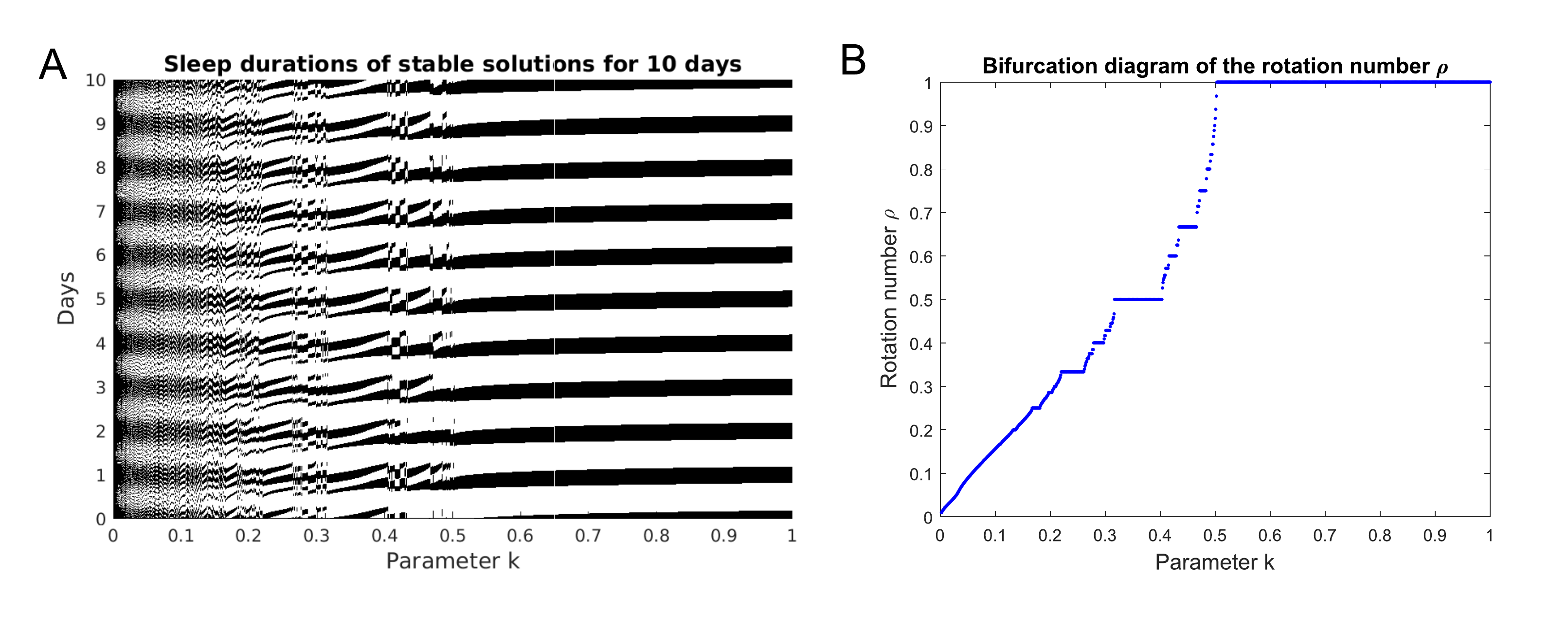}
       \caption{Multiple sleep episodes per day occur as time constants for the homeostatic sleep drive are decreased.  A: Simulated sleep periods (dark intervals) over the course of 10 days as the homeostatic sleep drive time constants are decreased by the scaling parameter $k$ (x-axis).  B: Bifurcation diagram of stable solutions in terms of the rotation number $\rho$ for the default parameter set. The parameter $k$ is on the $x-$axis and the rotation number $\rho$, defined as the number of circadian days over the number of sleep episodes in the stable sleep pattern is on the y-axis. The step size for $k$ was 0.001.}
        \label{fig:bfdiagrams}
    \end{figure}

     To examine how decreasing the homeostatic time constants $\tau_{hs}$ and $\tau_{hw}$ affects model solutions,  we introduce a scaling constant $k \in (0, 1]$ that multiplies both $\tau_{hs}$ and $\tau_{hw}$ in Eq. \ref{FF_h}. This is a simple scaling that preserves the ratio between time constants and is consistent with approaches in previous work \cite{skeldonphilips, skeldongaps}. We numerically computed model solutions (Figure \ref{fig:bfdiagrams}A) with respect to the bifurcation parameter $k$ to understand the change in the types of stable phase-locked solutions obtained as we decrease $k$ from 1. Specifically, we tracked the timing and duration of the sleep episodes of the stable solutions (black intervals) over the course of 10 days.
    
   We found that the number of sleep episodes per day increased as the time constants for the homeostatic drive decreased. At the default value $k=1$, the model produces one sleep episode per day (which we define here as a 24-hour cycle measured between two minimums of the circadian variable (c(t)). As we reduce $k$, the stability of solutions with one sleep episode per day is lost, and higher order patterns in which some days contain two sleep episodes may occur. In the interval ($k\in[0.317,0.403]$) a stable solution with two sleep episodes per day emerges. For smaller values of $k$, higher order patterns in which some days contain three sleep episodes appear, eventually resulting in the stable solution with three sleep episodes per day, and so on.     
 
      To quantify the sleep patterns associated with the attracting periodic orbits obtained for each value of $k$,  we define the rotation number, $\rho$, to be the number of circadian days $q$ over the number of sleep episodes $p$ occurring in one period of the stable orbit, i.e. $\rho=\frac{q}{p}$. Tracking $\rho$ as $k$ is decreased from 1 to 0 (Figure \ref{fig:bfdiagrams}B), we find that the rotation numbers vary as dictated by a Farey sequence \cite{granados,AvrutinBook}. In between neighboring intervals of $k$ displaying solutions with rotation numbers $\rho_1=\frac{a}{b}$ and $\rho_2=\frac{c}{d}$, where the greatest common divisor $\gcd(a,b)=1$,  $\gcd(c,d)=1 $, and $\lvert{ad-bc}\rvert=1$ is a $k$ interval with rotation number $\rho=\frac{a+c}{b+d}$. Such a Farey sequence of rotation numbers will generate a Devil's staircase-like structure for the rotation number as a function of $k$. This suggests that the stable solutions follow a period-adding bifurcation sequence that is consistent with previous work on systems governed by monotonic circle maps with discontinuities \cite{granados, skeldongaps, skeldon_dirks2019}.  Here, we numerically detect a subset of a Farey sequence of rotation numbers.

Interestingly, for small values of $k$ ($0 < k \le 0.18$), we obtain a denser set of rotation numbers from the computed solutions compared to the rotation numbers observed for larger values of $k$. For small $k$, the numerical results suggest that solutions with rotation numbers for all rational numbers less than about $\frac{1}{4}$ may exist. This is expected as sufficiently fast homeostatic time constants will result in continuous sleep onset circle maps. This occurs because fast time constants will prevent model trajectories from making tangent intersections with the saddle-node curves of the $Z-$shaped surface. In this case, the theory for monotonic, continuous circle maps guarantees that solutions exist with rotation numbers for all rational numbers \cite{granados,AvrutinBook}.

\subsubsection{Bifurcation sequences for emergence of stable solutions}
To identify the types of bifurcations leading to the gain (or loss) of stability of different sleep patterns for decreasing $k$, we track how the stable model trajectories and sleep onset maps evolve as we reduce $k$ for representative solutions associated with $\rho=1, \frac{2}{3}$, and $\frac{1}{2}$. Our analysis suggests that other stable solutions with $\rho \in [\frac{1}{2}, 1]$ will show the same bifurcation sequences. Solutions with smaller $\rho$ values may show different bifurcation sequences and are considered in Section \ref{dtocmap}.
    
    For $k=1$, the associated solution has rotation number $\rho=1$. As $k$ decreases, we describe the bifurcations associated with the loss of the $\rho=1$ solution. Similarly, we identify $k-$intervals associated with the existence of solutions with $\rho=\frac{2}{3}$ and $\frac{1}{2}$ and observe the bifurcation sequences associated with the emergence and loss of these solutions. These bifurcation sequences will include saddle-node (SN) bifurcations and border collision bifurcations of both stable (BC-S) and unstable (BC-U) fixed points of the maps. In the listing of the sequences, for all cases except $\rho=1$, the leftmost and rightmost bifurcations create and destroy, respectively, the stable solution with rotation number $\rho=\frac{q}{p}$ as $k$ is decreased.
    
    \subsubsection*{Border collision $\rightarrow$ saddle-node}  We first consider the loss of stability of the $\rho=1$ solution as $k$ is decreased from 1. The smallest value of $k$ for which this solution is stable is $k=0.503$ (Figure \ref{fig:BordersOfRhos}A,B). As $k$ is decreased towards this value, the stable periodic orbit shifts on the $Z-$shaped surface such that sleep onset occurs at earlier phases. The sleep onset map for $k=0.503$ reveals a saddle-node bifurcation to the right of the discontinuity (Figure \ref{fig:BordersOfRhos}A). The unstable fixed point associated with the saddle-node bifurcation was created at a higher value of $k$ ($k=0.504$) in a border collision bifurcation on the right side of the discontinuity (referred to as a Type I border collision in \cite{skeldongaps, skeldon_dirks2019}). Numerical simulations suggest that at the border collision bifurcation the unstable orbit makes a tangent intersection with the curve of saddle-node points on the upper plane of the $Z-$shaped surface. Thus, as $k$ decreases, the $\rho=1$ solution loses stability in the bifurcation sequence of
$$\textrm{BC-U} \rightarrow \textrm{SN}.$$ 

\subsubsection*{Saddle-node $\rightarrow$ border collision $\rightarrow$ border collision}
\begin{figure}[H]
\centering
        \includegraphics[scale=0.47]{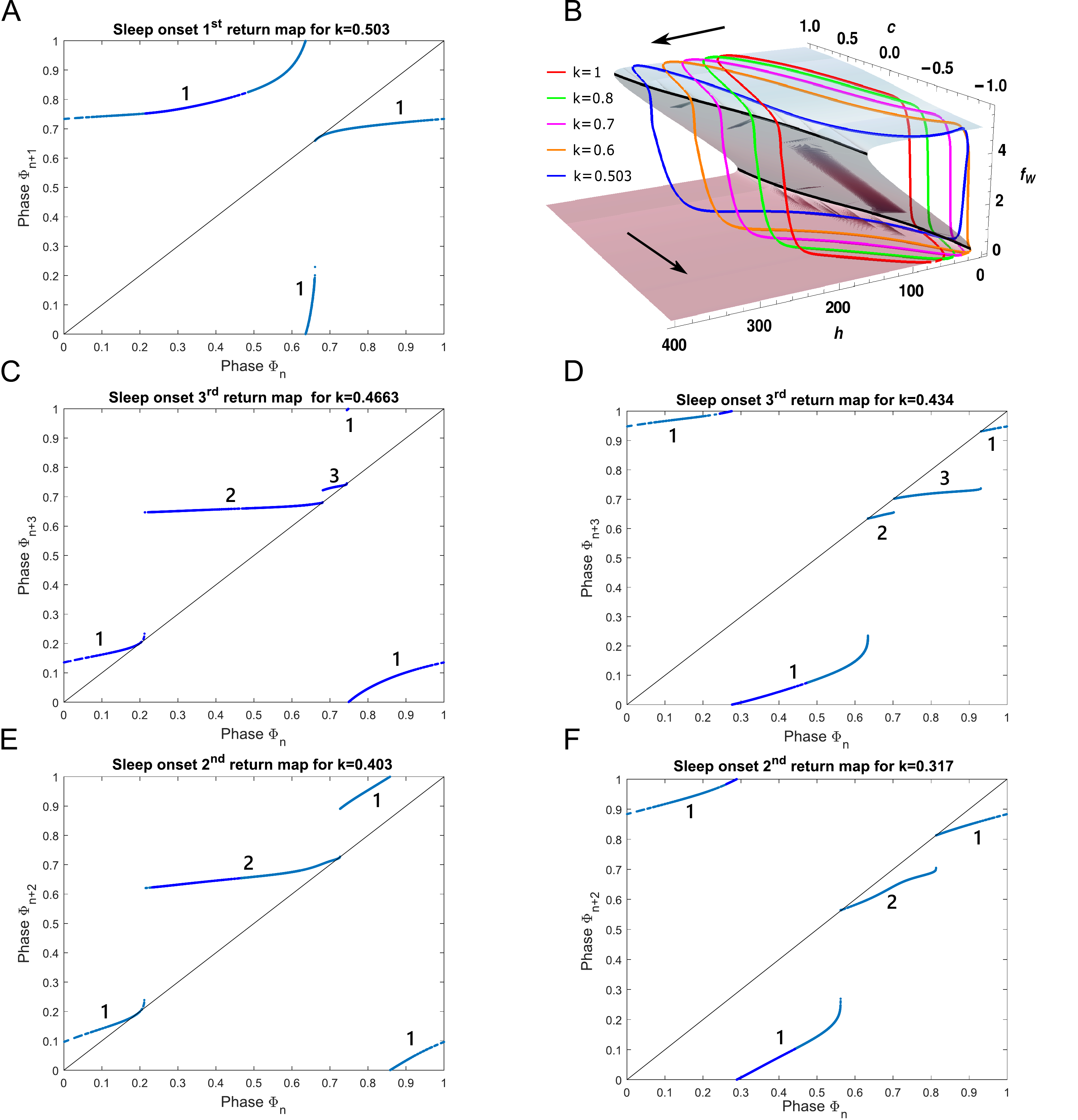}
    \caption{Sleep onset circle maps reveal the types of bifurcations at the emergence and disappearance of stable solutions as $k$ is varied. Distinct branches of the circle maps are labelled by numbers 1,2 and 3 as needed.
 A: The first return sleep onset map for $k=0.503$, the smallest value where the one sleep episode per day solution ($\rho=1$) is stable, shows a saddle-node bifurcation. 
 B: Evolution of stable $\rho=1$ periodic orbits plotted in relation to the $Z-$shaped surface in the $c-h-f_W$ space as $k$ approaches $k=0.503$. Each closed orbit corresponds to the stable solution for a particular value of the parameter $k$: $k=1$ (red), $k=0.8$ (green), $k=0.7$ (magenta), $k=0.6$ (orange), $k=0.503$ (blue). 
 C,D: Third return sleep onset maps for $k=0.4663$ (C) and $k=0.434$ (D). For this range of $k$ values the stable solution alternates between one and two sleep episodes per circadian cycle ($\rho = \frac{2}{3}$). The map has three branches (modulo 1) with a saddle-node bifurcation occurring at the right branch end at $k=0.4663$ (C) and a border collision occurring at the left branch end at $k=0.434$ (D). 
 E,F: The second return sleep onset maps for $k=0.403$ (E) and $k=0.317$ (F) between which exists the stable solution with two sleep episodes per circadian cycle. The map has two branches (modulo 1) with a saddle-node bifurcation occurring at the right branch end at $k=0.403$ (E) and a border collision occurring at the left branch end at $k=0.317$ (F). }
    \label{fig:BordersOfRhos}
\end{figure}

Next we describe the bifurcation sequences associated with the emergence and loss of a stable solution with alternating 1 and 2 sleeps per 24-h circadian cycle ($\rho = \frac{2}{3}$). This solution gains stability at $k=0.4663$ and loses stability at $k=0.434$. Fixed points associated with this solution appear in the 3rd return sleep onset map (Figure \ref{fig:BordersOfRhos}C,D). These maps consist of 3 separate branches, each showing an infinite slope at its right end and a finite slope at its left end (see Appendix \ref{Structure_of_maps}). Note that the two segments for lower $\Phi_{n+3}$ values form one connected branch (modulo 1) due to periodicity of the circle map. At $k=0.4663$, the map shows a saddle-node bifurcation near the infinite slope end of the map branches (numbered 1-3 in Figure \ref{fig:BordersOfRhos}C). The unstable fixed points are destroyed in a border collision (referred to as a Type II border collision in \cite{skeldongaps, skeldon_dirks2019}) at a slightly lower value of $k$ ($k=0.466$). Numerical solutions suggest that this border collision is associated with the unstable orbit making a tangent intersection with the upper curve of saddle-node points of the $Z-$shaped surface. As $k$ decreases to $k=0.434$, the map transitions so the stable fixed points move towards the left end of the map branches and disappear in a border collision bifurcation (Figure \ref{fig:BordersOfRhos}D). Numerical simulations indicate that this border collision bifurcation occurs due to a tangent intersection with the upper curve of saddle-node points of the $Z-$shaped surface. Thus, as $k$ is decreased, the emergence and disappearance of the $\rho=\frac{2}{3}$ stable solution occurs in the bifurcation sequence  $$\textrm{SN} \rightarrow \textrm{BC-U} \rightarrow \textrm{BC-S}.$$ 
    
We find that other stable solutions for lower values of $k$ also emerge through this same bifurcation sequence. For example, the $\rho=\frac{1}{2}$ solution with two sleep episodes per circadian cycle is stable in the interval $k \in [0.317,0.403]$. Fixed points for this solution appear in the second return sleep onset maps which in this regime consist of two separate branches, each with an infinite slope at the right end and a finite slope at the left end (Figure \ref{fig:BordersOfRhos}E,F). Again, as $k$ decreases the solution gains stability in a saddle-node bifurcation at the right end of the map branches where the unstable fixed points are destroyed in a border collision at $k=0.401$. Numerical simulations suggest that the unstable orbit makes a tangent intersection with the upper curve of saddle-node points of the $Z-$shaped surface. The fixed points disappear in a border collision with the left end of the map branches (Figure \ref{fig:BordersOfRhos}F), where again the stable orbit makes a tangent intersection with the upper curve of saddle-node points of the $Z-$shaped surface.


We have demonstrated that as $k$ varies, the appearance of tangent intersections of model trajectories with the curves of saddle-node points of the $Z-$shaped surface influence the occurrence of the bifurcations, and thereby, the emergence of stable orbits. The presence of tangent intersections depends, in part, on the circadian waveform and highlights the importance of the circadian drive on the bifurcation sequences.

\subsection{Varying the circadian waveform}
The circadian waveform reflects the time-varying profile of the firing rate of the SCN population. The properties of the SCN waveform are determined by interindividual differences as well as environmental light schedules that change with the seasons \cite{VANDERLEEST}. 
To investigate the effect of this waveform on the stable sleep-wake patterns, we varied the firing rate profile of the SCN population and tracked the existence of tangent intersections between model trajectories and the curves of saddle-node points of the $Z-$shaped surface. 
 
\begin{figure}[H]
		\centering
		\includegraphics[scale=0.39]{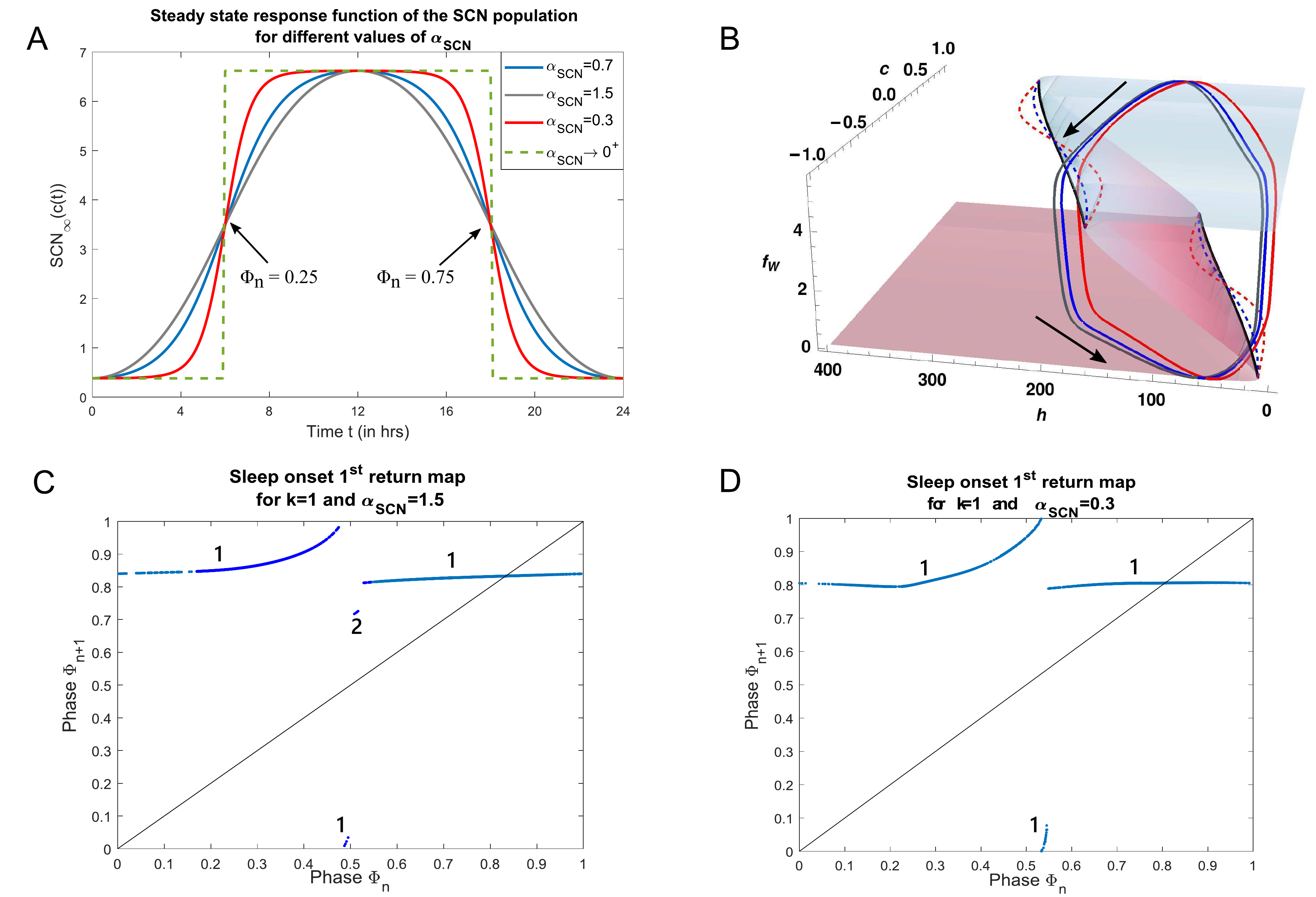}
	\caption{Effect of the parameter $\alpha_{SCN}$ on the circadian waveform, fast-slow decomposition surfaces and first return circle maps. A: Profile of $SCN_{\infty}(c(t))$ over 24 h for $\alpha_{SCN}$=0.7 (default value, blue), $\alpha_{SCN}$=1.5 (gray), $\alpha_{SCN}$=0.3 (red) and the limiting case $\alpha_{SCN}\to 0^{+}$ (dashed green). B: The $Z-$shaped surface of steady state solutions of the model fast subsystem showing the variation in the curve of saddle-node points with $\alpha_{SCN}$ ($\alpha_{SCN}=1.5$ (black), 0.7 (dashed blue) and 0.3 (dashed red)) with stable trajectories for $k=1$ ($\alpha_{SCN}=1.5$ (gray), 0.7 (blue) and 0.3 (red)).  C,D: First return sleep onset circle maps for $k=1$ and $\alpha_{SCN}=1.5$(C) and $0.3$(D). Distinct branches of the circle maps are labelled by the number 1 and 2 as needed.}
	\label{fig:aSCNeffect}
\end{figure}

Specifically, we modulated the circadian waveform such that its ``steepness'' (the transition region between low and high values of the SCN firing rate) varies without affecting the amplitude of the waveform. This is achieved by allowing the parameter $\alpha_{SCN}$ in the steady state response function of the SCN firing rate (Eq. \ref{FF_fSCN}) to vary from its default value $\alpha_{SCN}=0.7$ (Table \ref{tab:DefaultParameterSetTable}). We consider $\alpha_{SCN} \in (0,3]$. Decreasing or increasing $\alpha_{SCN}$ results in longer or shorter intervals, respectively, of high SCN firing rate activity (Figure \ref{fig:aSCNeffect}A) consistent with the response of SCN activity to longer or shorter environmental light periods \cite{Mrugala,VANDERLEEST}. 
  
To illustrate the effects of changing the profile of the $SCN_{\infty}(c(t))$ function on solutions with $k=1$ , we consider sleep onset maps and fast-slow decomposition for representative $\alpha_{SCN}$ values, $\alpha_{SCN}= 0.3$ and $1.5$ (Figure \ref{fig:aSCNeffect}B-D).
In the $Z-$shaped surface, the curves of saddle-node points have smaller (larger) curvature for larger (smaller) values of $\alpha_{SCN}$. Stable trajectories trace out similar hysteresis loops over the $Z-$surface.
In addition, the sleep onset first return maps have the same general shape as the default case, displaying a similar discontinuity with an infinite slope on its left side. The fixed points also occur at similar phases for $\alpha_{SCN}$ equal to 0.3, 0.7 and 1.5: namely $(\Phi_n, \Phi_{n+1}) = (0.8057, 0.8057)$, $(0.8242, 0.8242)$, and $(0.833,0.833)$, respectively. 

The map for $\alpha_{SCN}=1.5$ exhibits a second discontinuity resulting in a small map branch near $\Phi_n = 0.5$ (Figure \ref{fig:aSCNeffect}C). This discontinuity is caused by tangent intersections for trajectories associated with $\Phi_n$ values near 0.5: one initial phase produces a trajectory that makes a tangent intersection with the lower saddle-node curve, and a slightly higher initial phase produces the trajectory that makes a tangent intersection with the upper saddle-node curve. 
\subsection{Varying both homeostatic time constants and circadian waveform}


We study the combined effect of the parameter $\alpha_{SCN}\in(0,3]$ on the stable sleep-wake patterns obtained and bifurcation sequences arising as the homeostatic sleep drive time constants are scaled by $k$.  
To that end, we first consider the stable, phase-locked solutions obtained as $k$ is decreased for representative $\alpha_{SCN}$ values greater ($\alpha_{SCN} = 1.5$) and less ($\alpha_{SCN} = 0.3$) than the default value ($\alpha_{SCN} = 0.7$). We initially analyze the $\alpha_{SCN}$ effect on the bifurcation sequence for the loss of stability of the $\rho=1$ solution. Next, we compute the $(k,\alpha_{SCN})$ two-parameter bifurcation diagram to illustrate the evolution of bifurcation sequences over ranges of $k$ and $\alpha_{SCN}$ values. 

As we describe below, for lower values of $\alpha_{SCN}$, numerical simulations detect many fewer stable solutions associated with rotation numbers $\rho \in [\frac{1}{2}, 1]$. To verify this trend for the lowest values of $\alpha_{SCN}$, we additionally consider the limiting case of $\alpha_{SCN} \rightarrow 0^+$ corresponding to the SCN firing rate changing as a step function (see Section \ref{HSmodel}).

\subsubsection{Stable solutions for shallow and steep circadian waveforms (i.e., $\alpha_{SCN}=$1.5 and 0.3)}

\begin{figure}[H]
	\centering
 \includegraphics[scale=0.55]{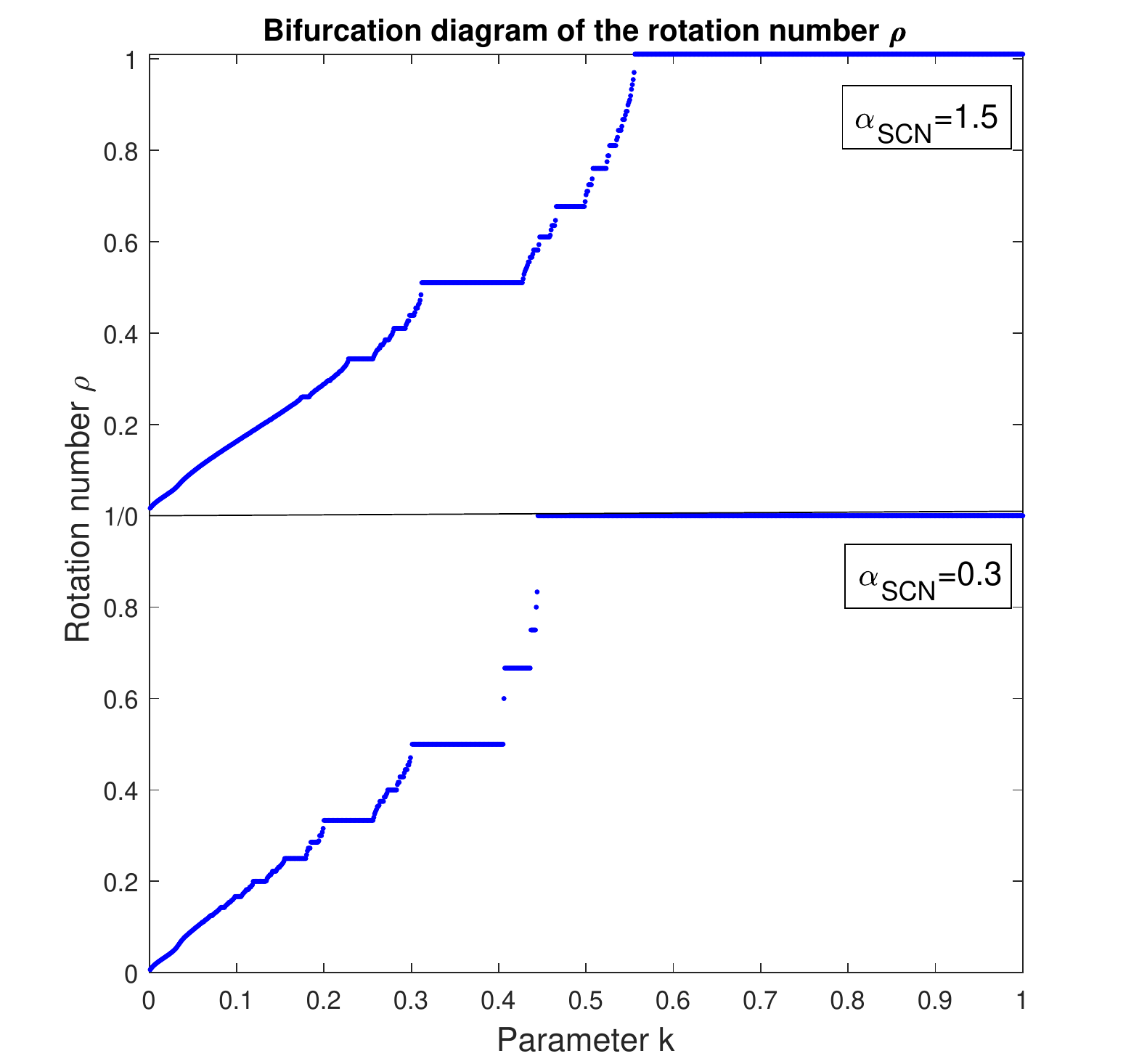}
 \caption{Comparison of the bifurcation diagrams of the rotation number $\rho$ for $\alpha_{SCN}=1.5$ (top) and $\alpha_{SCN}=0.3$ (bottom). Using a numerical approach to construct these diagrams, we obtained more types of periodic solutions with $\rho \in [\frac{1}{2},1]$ for larger $\alpha_{SCN}$ (shallower circadian waveform) compared to the solutions for smaller $\alpha_{SCN}$ (steeper circadian waveform). }
	\label{fig:ffdevilsstaicasea03and15}
\end{figure}

One key effect of changing the circadian waveform is that as $k$ is decreased from 1,  the $\rho=1$ solution corresponding to one sleep episode per circadian cycle loses stability earlier for larger values of $\alpha_{SCN}$ (shallower waveforms)   (Figure \ref{fig:ffdevilsstaicasea03and15}). For example, the $\rho=1$ solution loses stability at $k=0.556$, $k=0.503$, and $k=0.455$ for $\alpha_{SCN}=1.5, 0.7$ and $0.3$, respectively. As discussed above, the creation of tangencies of trajectories with the upper saddle-node curves of the $Z-$shaped surface is important in order for bifurcations to occur. Both parameters $k$ and $\alpha_{SCN}$ influence the creation of such tangent trajectories since the latter dictates the shape of the saddle-node curve and together they determine the angle at which a trajectory approaches it. 
As the upper saddle-node curve becomes steeper (for lower values of $\alpha_{SCN}$), $h$ must evolve faster for a trajectory orbit (stable or unstable) to become tangent to it, thus leading to the lower $k$ values when the $\rho=1$ solution loses stability.

Additionally, for larger values of $\alpha_{SCN}$, numerical simulations detect more stable solutions (than in the default $\alpha_{SCN}=0.7$ case) corresponding to distinct values of the rotation number, particularly types of $\rho=\frac{q}{p}$ periodic solutions within the intervals between $\rho=\frac{1}{p}$ periodic solutions (Figures \ref{fig:ffdevilsstaicasea03and15} (top) and \ref{fig:bfdiagrams}B). Conversely, for smaller $\alpha_{SCN}$ values, a winnowing (i.e., shrinking of the $k$-distance) of stable solutions with $\rho \in [\frac{1}{2}, 1]$ is observed (Figure \ref{fig:ffdevilsstaicasea03and15} (bottom)). While the arithmetic precision and the step size of the parameter $k$ in our numerical simulations could account for the inability to detect more solutions, we can conclude that stable solutions in this $\rho$ range exist over shorter $k$ intervals. 



\subsubsection{Bifurcations sequences for $\rho=1$ solutions for representative $\alpha_{SCN}$ values}

\subsubsection*{Border collision $\rightarrow$ saddle-node}

To understand the bifurcation leading to the loss of stability of the $\rho=1$ solution as $k$ decreases when $\alpha_{SCN}$=1.5, Figure \ref{fig:alldataa15and03}A displays the evolution of the stable periodic orbits for various values of $k$ ranging from $k=1$ to $k=0.556$, the $k$ value just before the loss of stability of the $\rho=1$ solution. As suggested by the absence of a tangent intersection of the trajectory with the saddle-node curve, the sleep-onset map demonstrates a saddle-node bifurcation near the right side of the discontinuity for the bifurcation value of $k$=0.556 (Figure \ref{fig:alldataa15and03}B). 
The unstable fixed point associated with the saddle-node bifurcation was created at a higher value of $k$ ($k=0.56$) in a border collision bifurcation on the map branch on the right side of the discontinuity. The associated unstable orbit makes a tangent intersection with the upper curve of saddle-node points of the $Z-$shaped surface. Similarly to the default $\alpha_{SCN}=0.7$ case, the $\rho=1$ solution for $\alpha_{SCN}>0.7$ loses stability in the bifurcation sequence of BC-U $\rightarrow$ SN as $k$ is decreased. 

\subsubsection*{Border collision} 

\begin{figure}[H]
\centering
    \includegraphics[scale=0.39]{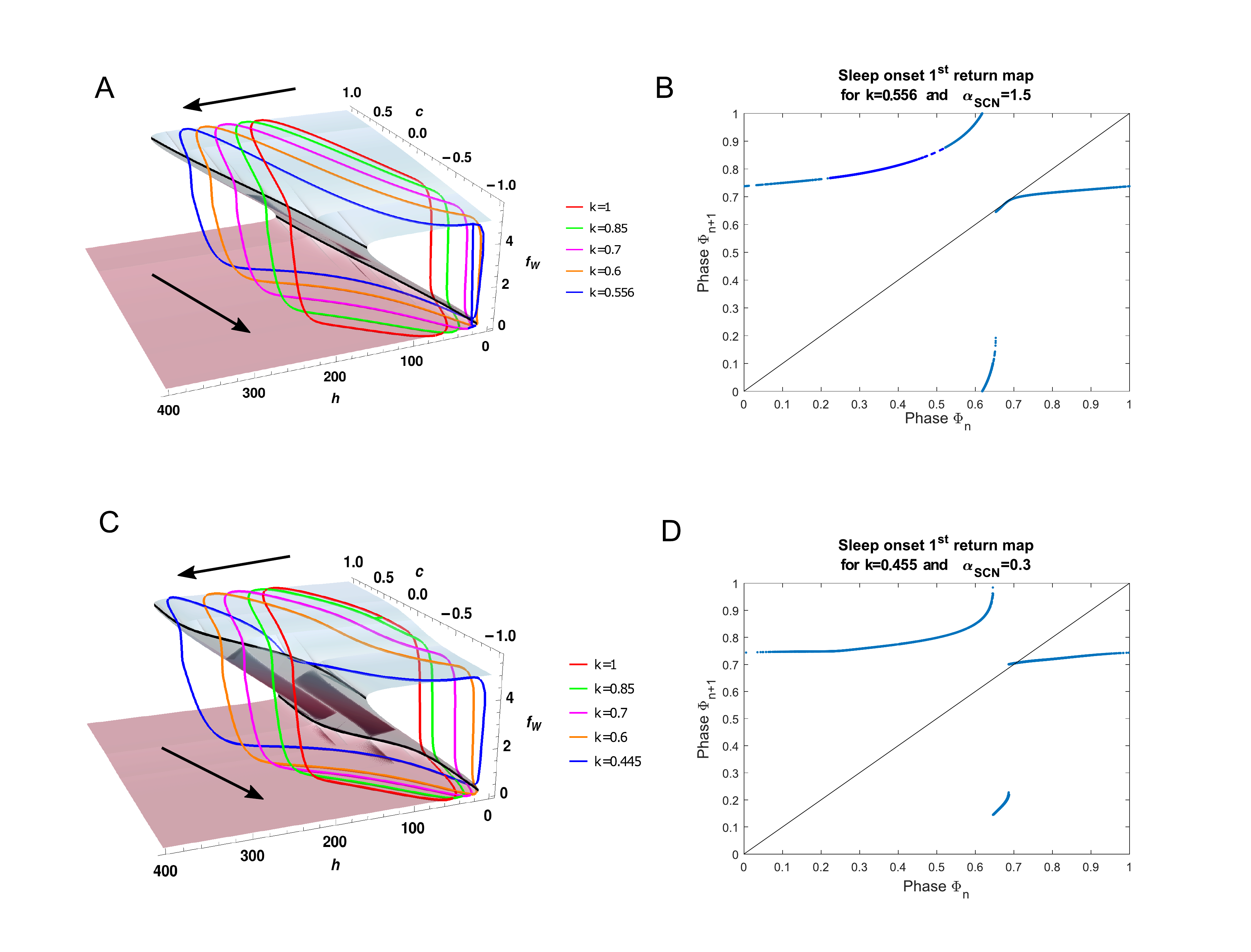}
	\caption{Bifurcations at the loss of stability of the $\rho=1$ solution for representative large (A,B) and small (C,D) $\alpha_{SCN}$ values. A: Stable trajectories for $\alpha_{SCN}=1.5$ and $k=1$ (red), $k=0.8$ (green), $k=0.6$ (orange) and $k=0.556$ (blue). B: First return sleep onset map for $\alpha_{SCN}=1.5$ and $k=0.556$ indicates the loss of stability of the $\rho=1$ solution occurs due to a saddle-node bifurcation. C: Stable trajectories for $\alpha_{SCN}=0.3$ and $k=1$ (red), $k=0.8$ (green), $k=0.6$ (magenta), $k=0.55$ (orange) and $k=0.506$ (blue). D: First return map for $\alpha_{SCN}=0.3$ and $k=0.445$ indicates the $\rho=1$ solution loses stability due to a border collision bifurcation.}
    \label{fig:alldataa15and03}
\end{figure}

For $\alpha_{SCN}=0.3$, we observe a different bifurcation sequence when the $\rho=1$ solution loses stability at $k=0.445$. At this bifurcation point, the  associated sleep-onset map continues to show two discontinuities as observed for the map for $k=1$.  The map demonstrates that the bifurcation occurs due to a border collision on the right side of the discontinuity (Figure \ref{fig:alldataa15and03}D). This border collision corresponds to the stable trajectory (Figure \ref{fig:alldataa15and03}C, blue curve) creating a tangency at the upper saddle-node curve of the $Z-$shaped surface (Figure \ref{fig:alldataa15and03}C). This suggests that for smaller $\alpha_{SCN}$ values, the $\rho=1$ solution ceases to exist due to a BC-S bifurcation. 

\subsubsection{Bifurcation sequences in ($k,\alpha_{SCN})$ parameter space}

\begin{figure}[H]
	\centering
	\includegraphics[scale=0.575]{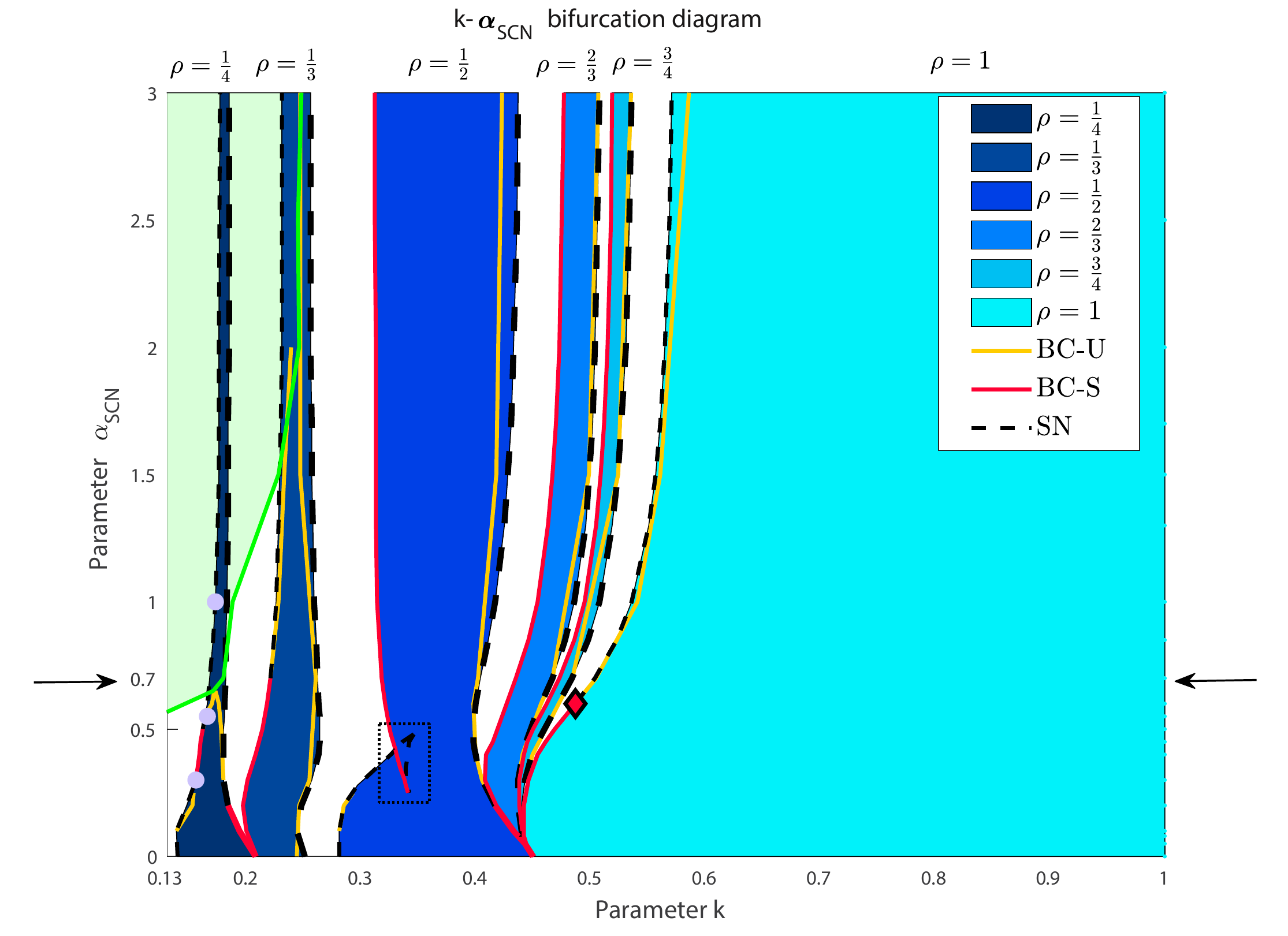}
	\caption{Two parameter bifurcation diagram with respect to $k$ and $\alpha_{SCN}$. Colored areas indicate parameter regions (or tongues) where the following stable, phase-locked solutions exist (from left to right): $\rho=\frac{1}{4}, \frac{1}{3}, \frac{1}{2}, \frac{2}{3}, \frac{3}{4} ,1$. Line type indicates bifurcation type: saddle-node (dashed black), border collision of a stable fixed point (solid red), border collision of an unstable fixed point (solid yellow). Diamond indicates transition between bifurcation sequences governing loss of stability of the $\rho=1$ solution (see Figure \ref{fig:alldataa15and03}). Arrows indicate the default $\alpha_{SCN}$ value of 0.7. The green line is the set of $(k,\alpha_{SCN})$ points that forms the boundary between regions where maps are discontinuous and continuous, and the light green shaded region indicates the $(k, \alpha_{SCN})$ values for which the map is continuous. The black dotted rectangle indicates a zoomed in version of the two-parameter bifurcation diagram shown in Figure \ref{fig:RhoOp5_Bistability}. The three light purple bullets indicate the $(k,\alpha_{SCN})$ values of the maps shown in Figure \ref{fig:1over4tongue}.}
	\label{fig:kalphaSCNbifdiagram}
\end{figure}

To illustrate the evolution of bifurcation sequences over a range of homeostatic time constants and circadian waveforms, we constructed a two-parameter bifurcation diagram with respect to $k$ and $\alpha_{SCN}$ (Figure \ref{fig:kalphaSCNbifdiagram}). The $\rho=1$ entrainment region (cyan) is bordered on the left by a curve of $(k,\alpha_{SCN})$ values associated with stable fixed points at which a saddle-node bifurcation (dashed black) or a border collision (solid red) occurs. 

The transition from the regime where the stable $\rho=1$ solution is lost due to the BC-U $\rightarrow$ SN bifurcation sequence to the regime where it is lost due to a BC-S bifurcation occurs continuously as $\alpha_{SCN}$ decreases. In the  BC-U$\rightarrow$SN regime, at the $k$ value associated with the border collision the slope of the map at the created fixed point is greater than 1 (but finite) resulting in an unstable fixed point. As $\alpha_{SCN}$ decreases, the slope of the map curve at the unstable fixed point created in this bifurcation also decreases. 

The two regimes are separated at $(k,\alpha_{SCN}) = (0.486,0.6)$ which is marked with a diamond. At this point, the curve of stable fixed points (solid red)  merges with the curve of unstable fixed points (solid yellow). For $\alpha_{SCN}>0.6$, the unstable fixed points are created in a border collision bifurcation occurring at a higher $k$ value than the $k$ value associated with the saddle-node bifurcation that forms the boundary of the $\rho=1$ entrainment region. The transition between bifurcation regimes occurs at  $(k,\alpha_{SCN}) = (0.486,0.6)$. Here, the fixed point of the map coincides with the end point (border) of the map curve, and the slope of the map curve at that point is equal to 1. For $\alpha_{SCN}<0.6$, the stable fixed point associated with the $\rho=1$ solution is lost directly due to a border collision bifurcation. 

The two parameter bifurcation diagram of Figure \ref{fig:kalphaSCNbifdiagram} also shows the entrainment regions (or tongues) in $(k, \alpha_{SCN})$ space for stable solutions with $\rho=\frac{3}{4}, \frac{2}{3}$ and $\frac{1}{2}$. For $\alpha_{SCN} \ge 0.2$, each of these solutions gains stability, as $k$ is decreased, through a  saddle-node bifurcation (dashed black curve) that is followed by a border collision (solid yellow curve) that eliminates an unstable fixed point. For the $\rho=\frac{3}{4}, \frac{2}{3}$ solutions as $k$ is further decreased, stability is lost through a border collision (solid red curve) resulting in the bifurcation sequence $\textrm{SN}\rightarrow \textrm{BC-U} \rightarrow \textrm{BC-S},$ similar to the bifurcation sequence observed for these solutions for the default $\alpha_{SCN} = 0.7$.

While not computed explicitly, we argue that this bifurcation sequence delimits the entrainment regions of all solutions with $\rho \in (\frac{1}{2}, 1)$ and $\alpha_{SCN} \ge 0.2$.  Specifically, for a stable solution with  rotation number $\rho=\frac{q}{p}$ we consider the $p^{th}$ order return map. As explained in Appendix \ref{Structure_of_maps}, the $p^{th}$ order map retains similar structure as the first return map. For example, for values of $(k,\alpha_{SCN})$ where the first return map is discontinuous, the $p^{th}$ return map is likewise discontinuous. Furthermore, the $p^{th}$ return map has $p$ discontinuities corresponding to each discontinuity in the first return map. For $\alpha_{SCN} \ge 0.2$ and all the $k$ values where solutions with these $\rho$ values exist, first return maps display an infinite slope at the left of a discontinuity and a finite slope on the right. The $p^{th}$ return map similarly shows this structure in each of the branches of the map. Computing maps at the $k$ values where these solutions gain and lose stability reveals that stable fixed points are created on map branches to the left of a discontinuity and are lost on map branches to the right of a discontinuity.
Specifically, for fixed $\alpha_{SCN}$, we consider the $p^{th}$ return map at the highest value of $k$ for which $\rho=\frac{q}{p}$ exists. On the $p$ branches associated with this map, there are $p$ saddle-node points formed by the infinite slope end of the map branches. As $k$ is decreased, $p$ unstable fixed points are eliminated in a border collision bifurcation on the infinite slope ends of the $p$ map branches. As $k$ is decreased further, the $p$ stable fixed points for the $\rho=\frac{q}{p}$ solution disappear in a border collision bifurcation at the finite slope end of the $p$ map branches. Since, for decreasing $k$, the bifurcation sequence SN $\rightarrow$ BC-U $\rightarrow$ BC-S is predicted by the structure of the map, we expect that all solutions with $\rho \in (\frac{1}{2}, 1]$ and $\alpha_{SCN} \ge 0.2$ will show a similar bifurcation sequence.

\begin{figure}[hbt!]
    \centering
   \hspace{0cm} 
   \includegraphics[scale=0.45]{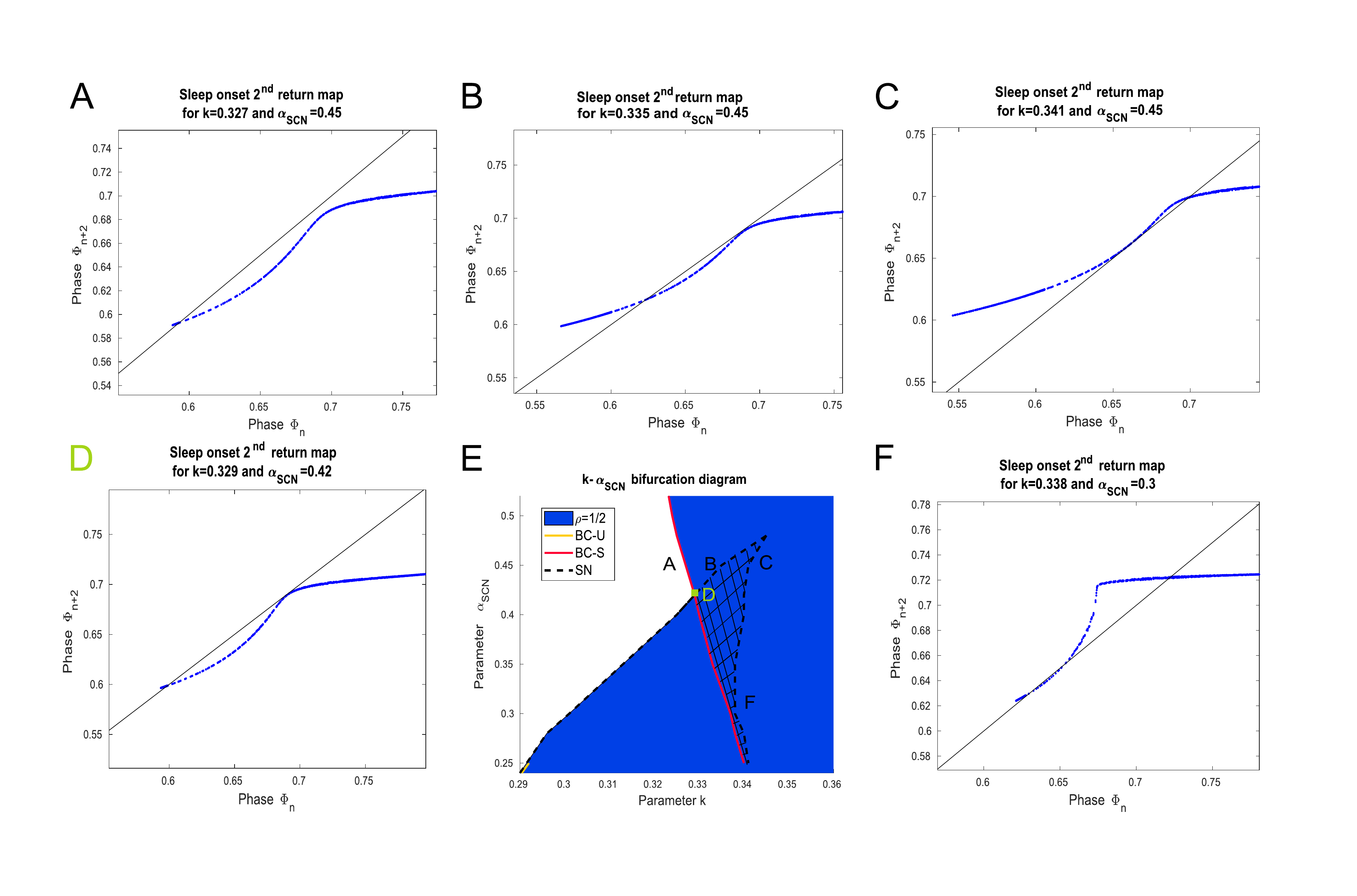}
    \caption{Bifurcations creating a bistability island in the $\rho=\frac{1}{2}$ stable entrainment region. A-C: Evolution of the second return map for $\alpha_{SCN}=0.45$. Here only one branch of the map is shown. For $k=0.341$ (C) a saddle-node bifurcation occurs at the lower part of the map curve. This gives birth to a new pair of stable and unstable fixed points. At $k=0.335$ (B) another saddle-node bifurcation leads to the collision of the new unstable and original stable fixed points. Complete loss of stability of the $\rho=\frac{1}{2}$ solution occurs at $k=0.327$ (A) in a border collision. D: At $\alpha_{SCN}=0.42$, the loss of bistability coincides with the loss of stability of the $\rho=\frac{1}{2}$ solution. At $k=0.329$ a saddle-node and a border collision eliminate two stable and one unstable fixed points. E: Close up of the two-parameter bifurcation diagram in $(k, \alpha_{SCN})$ space shown in Figure \ref{fig:kalphaSCNbifdiagram} shows the bistability island within the $\rho = \frac{1}{2}$ stable entrainment region. Letters in panel E correspond to maps shown in panels A-D and F. F: Second return map curve showing a ``sharp cornered S'' shape for $\alpha_{SCN}=0.3$ and $k=0.338$ where a saddle-node bifurcation initiates the interval of bistability.  
    }
      \label{fig:RhoOp5_Bistability}
\end{figure}

The bifurcations bounding the $\rho=\frac{1}{2}$ entrainment region are the same for $\alpha_{SCN} > 0.42$. However, the bifurcation governing the loss of stability of the $\rho=\frac{1}{2}$ solution changes to a saddle-node for $\alpha_{SCN} < 0.42$ (Figure \ref{fig:kalphaSCNbifdiagram}). This exchange in the bifurcations is a result of a small region or ``island'' of bistability emerging in the interior of the $\rho=\frac{1}{2}$ entrainment tongue (Figure \ref{fig:RhoOp5_Bistability}E). The ``bistability island'' exists for $\alpha_{SCN}\in[0.25, 0.48]$. It is bounded by curves of saddle-node bifurcations for high $k$ values while for lower $k$, it is bounded by a saddle-node curve for $\alpha_{SCN} \in (0.42, 0.48)$ and a curve of border collisions for $\alpha_{SCN} \in [0.25, 0.42)$. As described below, at $\alpha_{SCN}=0.42$, both bifurcations  occur at the same value of $k$, enabling the switch in bifurcation type eliminating the stable $\rho=\frac{1}{2}$ solution.  

The region of bistability occurs due to the curves of the second return map becoming S-shaped which allows for multiple intersections with the diagonal $\Phi_{n+2}=\Phi_{n}$, and thus multiple fixed points. Specifically, in this $\alpha_{SCN}$ interval, as $k$ decreases within the $\rho = \frac{1}{2}$ entrainment interval, the second pair of stable fixed points (and a pair of unstable fixed points) are created in the second return map due to a saddle-node bifurcation at the lower knees of the S-shaped map curves (Figure \ref{fig:RhoOp5_Bistability}C, figure shows one of the map branches). On the map branch shown in Figure \ref{fig:RhoOp5_Bistability}C, the original $\rho = \frac{1}{2}$ solution corresponds to the stable fixed point at higher sleep onset phase and the newly created solution with the stable fixed point at lower sleep onset phase. The newly created unstable solution has a sleep onset phase between those of the stable fixed points. For $\alpha_{SCN} \in (0.42,0.48]$, as $k$ decreases further, the new unstable fixed points and the original stable fixed points approach each other and eventually collide in a saddle-node bifurcation  at the upper knees of the S-shaped map curves (Figure \ref{fig:RhoOp5_Bistability}B).
This bifurcation marks the end of the interval of bistability and the newly created stable fixed points remain.  These fixed points are eliminated, and the $\rho = \frac{1}{2}$ solution loses stability, in a border collision at the left ends of the map branches (Figure \ref{fig:RhoOp5_Bistability}A). Thus, the complete bifurcation sequence for $\alpha_{SCN} \in (0.42,0.48]$ is $$\textrm{SN} \rightarrow \textrm{BC-U} \rightarrow \textrm{SN}\rightarrow \textrm{SN} \rightarrow \textrm{BC-S}.$$

At $\alpha_{SCN}=0.42$, the end of the bistability interval coincides with the loss of stability of the $\rho = \frac{1}{2}$ solution as the saddle-node bifurcation at the upper knees of the S-shaped map curves occurs at the same $k$ value as the border collision at the left ends of the map branches  (Figure \ref{fig:RhoOp5_Bistability}D, $k=0.329$). At this value of $k$, the two pairs of stable fixed points (one on each map branch of the second return map) lose stability simultaneously. The fixed points corresponding to the original stable solution that initiated the $\rho=\frac{1}{2}$ tongue (at higher sleep onset phase in Figure \ref{fig:RhoOp5_Bistability}D) lose stability due to a saddle-node bifurcation with the unstable fixed points. The other stable fixed points (at lower sleep onset phase in Figure \ref{fig:RhoOp5_Bistability}D)
lose stability due to a border collision. This causes the exchange of the bifurcation dictating the loss of stability of the $\rho=\frac{1}{2}$ solution from a border collision to a saddle node. Thus at $\alpha_{SCN}=0.42$ the full bifurcation sequence is $$\textrm{SN} \rightarrow \textrm{BC-U} \rightarrow \textrm{SN} \rightarrow \textrm{BC-S}+ \textrm{SN}.$$

For $\alpha_{SCN} \in [0.25,0.42)$, the shape of the map branches in the second return map transitions to a ``sharp cornered S'' (Figure \ref{fig:RhoOp5_Bistability}F). When the map is continuous in this sharp cornered S shape, the following bifurcation sequence takes place  $$\textrm{SN} \rightarrow \textrm{BC-U} \rightarrow \textrm{SN} \rightarrow \textrm{BC-S}\rightarrow \textrm{SN}.$$
Here, the stable fixed points that introduced bistability (at lower phase in the figure) lose stability first in a border collision at the left end of the map curves. For lower $k$ values, the unstable and original stable fixed points collide in a saddle-node bifurcation which eliminates the stable  $\rho=\frac{1}{2}$ solution.

As $\alpha_{SCN}$ approaches $0.25$, a discontinuity can occur in the map in this sharp cornered S shape, where the slope of the map branches are infinite from the left and finite from the right. In this case, the complete bifurcation sequence is $$\textrm{SN} \rightarrow \textrm{BC-U} \rightarrow \textrm{SN}\rightarrow\textrm{BC-S}\rightarrow \textrm{BC-U} \rightarrow\textrm{SN}.$$ The last border collision bifurcation creates another pair of unstable fixed points (one on each of the associated branches of the second return map), as the sharp cornered S shape of the map curves starts deforming as $k$ is decreased.

\subsection{Bifurcation sequences for small $k$}
\label{dtocmap}
 As noted previously, for small values of $k$, the homeostatic sleep drive varies more quickly relative to the SCN firing rate especially for high values of $\alpha_{SCN}$, thereby making tangent intersections of the solution trajectory with the curve of saddle-node points on the $Z-$surface less likely. As a result, the associated sleep onset maps can be continuous. This affects the bifurcation sequence delimiting stability of solutions with rotation numbers $\rho \le \frac{1}{3}$. The two parameter bifurcation diagram  can be separated into regimes associated with continuous or discontinuous sleep onset maps. There exists a curve of $(k,\alpha_{SCN})$ points (Figure \ref{fig:kalphaSCNbifdiagram}, solid green curve) above which the map is continuous (Figure \ref{fig:kalphaSCNbifdiagram}, light green area). We will refer to this (green) curve as the transition zone. 
 
 We note that bifurcation sequences closer to the continuous regime, and hence across the transition zone, may not involve border collision bifurcations associated with the creation or destruction of stable fixed points (BC-S). As the maps obtain larger discontinuities, we observe bifurcation sequences similar to the ones delimiting the entrainment regions we have encountered so far. We describe representative examples of the bifurcations across the transition zone with the $\rho=\frac{1}{4}$ solutions for different values of $\alpha_{SCN}$. 

\subsubsection*{Saddle-node $\rightarrow$ saddle-node} 
\begin{figure}[H]
     \centering
   \hspace{-0.75cm}
   \includegraphics[scale=0.385]{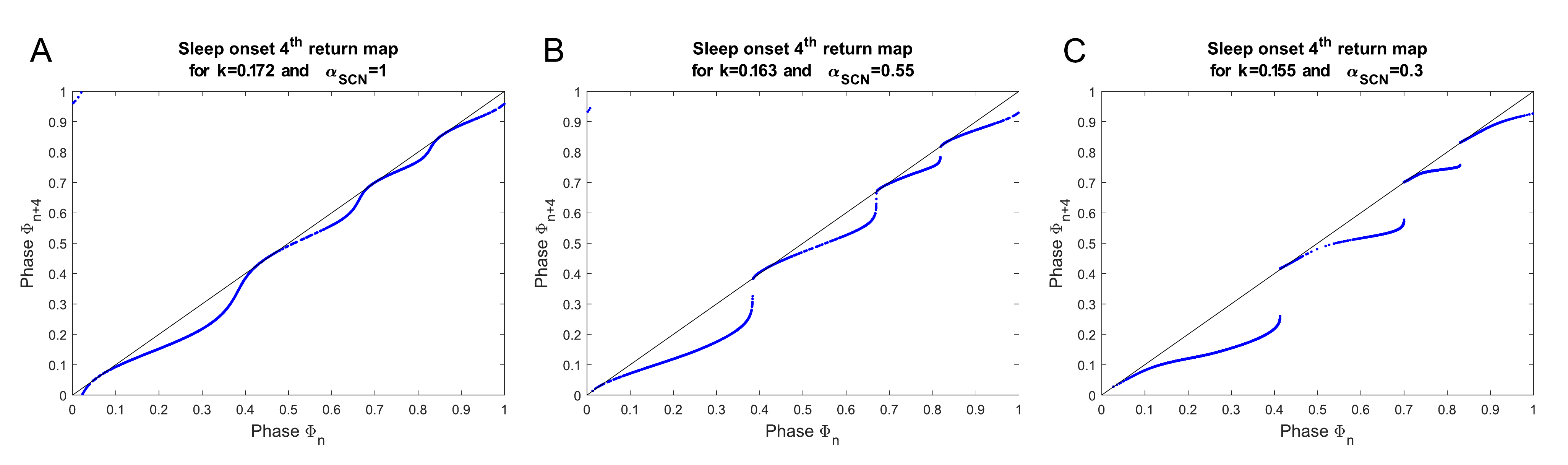}
     \caption{Transition from a continuous to a discontinuous sleep onset map within the $\rho=\frac{1}{4}$ entrainment region. Fourth return sleep onset maps are shown at smallest $k$ values where the $\rho=\frac{1}{4}$ solution is stable for different $\alpha_{SCN}$ values.  A: For $\alpha_{SCN}=1$ the map is continuous and loss of the stable $\rho=\frac{1}{4}$ solution is due to a saddle-node bifurcation. B: For $\alpha_{SCN}=0.55$ the map is discontinuous, but the slope of the map branches on the right of each discontinuity is greater than 1 in magnitude, leading to a border collision that generates an unstable fixed point followed by a saddle-node bifurcation as $k$ decreases. C: For $\alpha_{SCN}=0.3$ the map is discontinuous. A border collision on the right of each discontinuity leads to loss of the stable fixed points associated with the $\rho=\frac{1}{4}$ solution.}
     \label{fig:1over4tongue}
 \end{figure}

 For pairs of $(k,\alpha_{SCN})$ values above the transition zone (see light green area, the first and fourth return sleep onset maps are continuous, and hence saddle-node bifurcations lead to loss of stability of the $\rho=\frac{1}{4}$ periodic solution (Figure \ref{fig:1over4tongue}A). The fourth return map has four pairs of stable and unstable fixed points, and the unstable fixed points remain over the $k$ interval where the solution is stable, i.e. the unstable fixed point is not lost through a border collision. Thus, the stability of the $\rho=\frac{1}{4}$ periodic solution in this regime occurs in the bifurcation sequence $$\textrm{SN} \rightarrow \textrm{SN}.$$

\subsubsection*{Saddle-node $\rightarrow$ border collision $\rightarrow$ border collision [$\rightarrow$ saddle-node]} For ($k,\alpha_{SCN}$) pairs below the transition zone, the first and fourth return maps are discontinuous with infinite slopes to the left of the discontinuity and finite slopes to the right of the discontinuity. A saddle-node bifurcation leads to gain of stability of the $\rho=\frac{1}{4}$ periodic solution as $k$ is decreased. This is followed by a border collision bifurcation at a slightly lower value of $k$. 

For larger $\alpha_{SCN}$ values, right below the transition zone, the slope of the discontinuous map is greater than 1 at the right of the discontinuity (Figure \ref{fig:1over4tongue}B). Therefore, for decreasing $k$ there is first a border collision bifurcation that generates an unstable fixed point, and then this unstable fixed point eventually collides with the stable fixed point in a saddle-node bifurcation as $k$ decreases. Thus, stability of the $\rho=\frac{1}{4}$ solution in this regime occurs in the bifurcation sequence of $$\textrm{SN} \rightarrow \textrm{BC-U} \rightarrow \textrm{BC-U} \rightarrow  \textrm{SN}.$$ 
As $\alpha_{SCN}$ is reduced, the slope of the map decreases smoothly to values less than 1. Then the stable fixed point ceases to exist due to a border collision bifurcation as $k$ decreases. In particular, the full bifurcation sequence is
  $\textrm{SN} \rightarrow \textrm{BC-U} \rightarrow \textrm{BC-S}.$
  
   Transitions between these bifurcation sequences occurred smoothly. Specifically, for the $\rho=\frac{1}{4}$ solution at the transition zone, the unstable fixed points associated with the border collisions appear at the same $k$ value and their $k$ values diverge for smaller $\alpha_{SCN}$. Additionally, for the $\rho=\frac{1}{4}$ and $\frac{1}{3}$ solutions, as the bifurcation at the loss of stability as $k$ decreases changes, the $k$ values at the border collisions associated with the creation of unstable fixed points (solid yellow curve) merges with the values associated with the destruction of the stable fixed points (solid red curve) as $\alpha_{SCN}$ decreases, as we have observed with solutions of $\rho>\frac{1}{3}$. We expect that similar bifurcation sequences delimit other stable solutions with $\rho < \frac{1}{3}$ across the transition zone, as the associated maps are expected to maintain a similar structure. \textcolor{red}{}
 

\section{Circadian Hard-Switch model in limit $\alpha_{SCN}\to 0^{+}$}
\label{HSmodel}
As shown in the $k$-$\alpha_{SCN}$ bifurcation diagram (Figure \ref{fig:kalphaSCNbifdiagram}), as $\alpha_{SCN}$ decreases, 
the widths of the $k$ intervals associated with each stable entrainment regime change in the following way: For some rotation numbers, such as $\rho = 1, \frac{1}{2},$ and $\frac{1}{4}$, the $k$ intervals  expand for low $\alpha_{SCN}$ values, while for $\rho = \frac{3}{4}$, $\frac{2}{3}$ and $\frac{1}{3}$, the $k$ intervals contract. In the limit $\alpha_{SCN}\to 0^{+}$, this leads to the loss of stable solutions with $\rho \in (\frac{1}{2}, 1)$ and $\rho \in (\frac{1}{4}, \frac{1}{3})$, and a change in the bifurcation sequence bordering the $\rho = \frac{1}{2}$ and $\rho = \frac{1}{4}$ stable solutions. To analyze this change in the size of k-intervals for a small $\alpha_{SCN}$, we consider the model in the limit $\alpha_{SCN}\to 0^{+}$. We refer to this model as the circadian hard switch (CHS) model. 

In this section, we first formally define the CHS model and then describe the stable solutions obtained as $k$ decreases from 1 with a particular focus on the bifurcations delimiting the stable $\rho=\frac{1}{2}, \frac{1}{3}$ and $\frac{1}{4}$ solutions. Based on how the $\rho=1$ (and $\rho=\frac{1}{3}$) solutions directly transition to the $\rho=\frac{1}{2}$ (and $\rho=\frac{1}{4}$) solutions in the CHS model, allows us to explain why the $k$ intervals for stable solutions with $\rho \in (\frac{1}{2}, 1)$ (and $\rho \in (\frac{1}{4}, \frac{1}{3})$) shrink for small $\alpha_{SCN}$.

\subsection{Definition of the Filippov system with two switching boundaries}

In the limit as $\alpha_{SCN}\to 0^{+}$, the firing rate response function of the SCN population can be approximated by a step function. This introduces a second discontinuity in the $f_{SCN}$ derivative, when $c$ crosses $\beta_{SCN}$:

\begin{equation}\label{HSa_fSCN}
\dfrac{df_{SCN}}{dt}=\dfrac{SCN_{max}\cdot0.5\cdot\bigg(1+\tanh\Big(\frac{1}{0.7}\Big)(2\mathcal{H}(c-\beta_{SCN})-1)\bigg)-f_{SCN}}{\tau_{SCN}}\,,
\end{equation}

where $\mathcal{H}$ is the Heaviside function. Then, our model becomes a Filippov system with two switching boundaries \cite{FilippovBook}: one represents the switch between sleep and wake, and the other represents a switch between high and low activity in the circadian drive $c(t)$ as occurs in SCN firing rate over the 24 hour day \cite{VANDERLEEST}. 

\begin{figure}[H]
		\centering
	\includegraphics[scale=0.42]{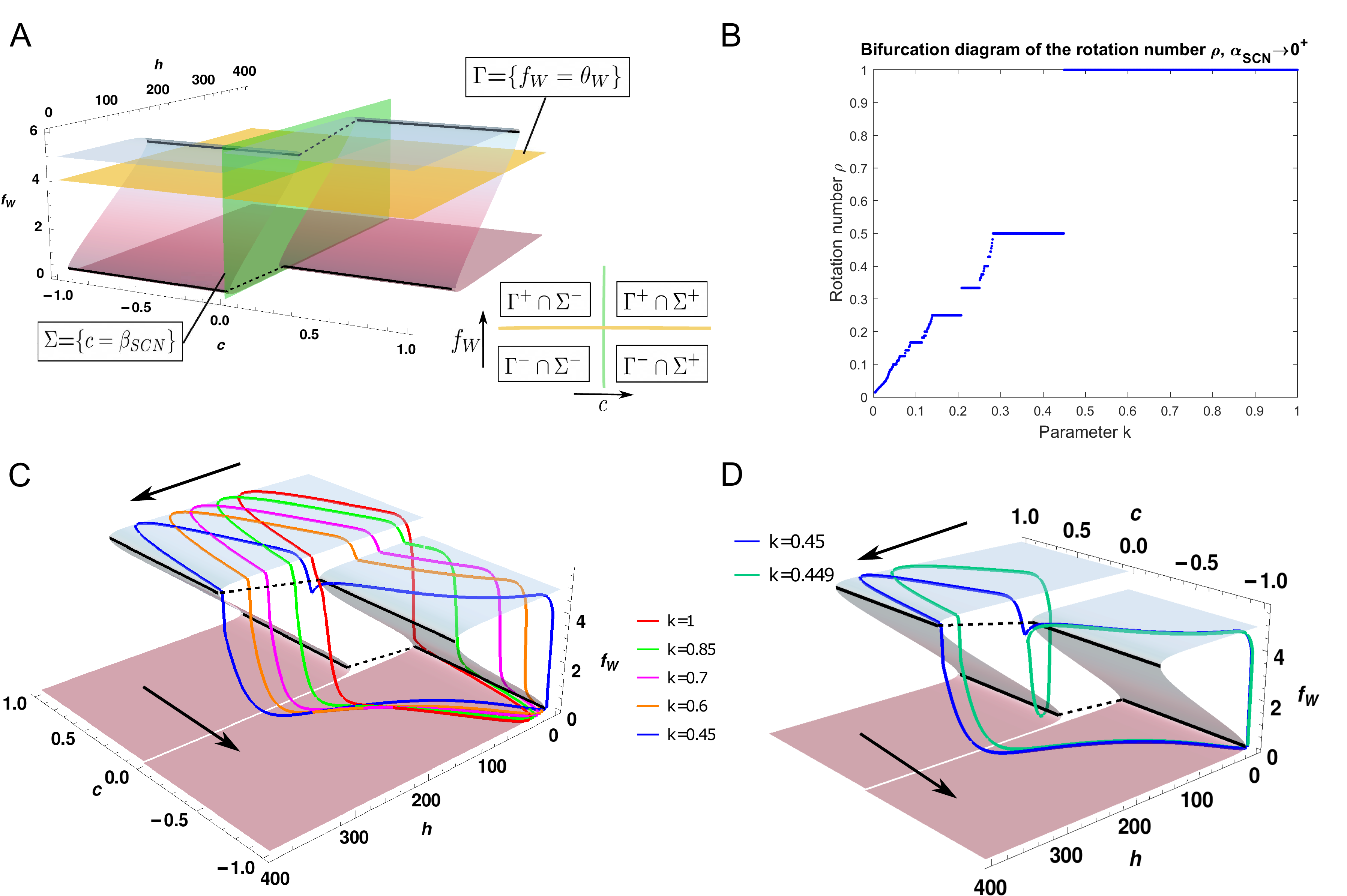}
\caption{Dynamics and bifurcation structure in the circadian hard switch (CHS) model. A: The CHS model is a Filippov system with two boundaries, $\Gamma$ (yellow plane) and $\Sigma$ (green plane). Therefore, in the $c-h-f_W$ space we can visualize the fast-slow surface being divided into four regions, $\Gamma^{+}\cap\Sigma^{+}$, $\Gamma^{+}\cap\Sigma^{-}$, $\Gamma^{-}\cap\Sigma^{-}$ and $\Gamma^{-}\cap\Sigma^{+}$. We have plotted the individual fast-slow surfaces for the dynamical system when $c>\beta_{SCN}$ (corresponds to $\Sigma^{+}$) and $c<\beta_{SCN}$ (corresponds to $\Sigma^{-}$). In each of these regions, the system is smooth, but a discontinuity in the derivative occurs as the system crosses a boundary. B: The bifurcation diagram of the rotation number $\rho$ for the CHS model. C: The evolution of the stable solutions leading to the loss of stability of the $\rho = 1$ solution. In this regime, a sleep onset always occurs at the same circadian phase when the trajectory crosses the boundary $\Sigma$. D: Decreasing the value of the scaling parameter $k$ leads to loss of stability of the $\rho = 1$ solution and emergence of the stable $\rho = \frac{1}{2}$ solution. }
\label{fig:HS_figures}
\end{figure}

To define the circadian hard switch (CHS) model, we introduce the new switching boundary $\Sigma$ in addition to the original switching boundary $\Gamma$ where $\Sigma$ is defined as:
$\Sigma=\{c=\beta_{SCN}\}$, where $\beta_{SCN}=0$ (Figure \ref{fig:HS_figures}A).

\noindent The regions lying on either side of each boundary are then defined as:
\begin{itemize}
    \item $\Sigma^{+}=\{c>\beta_{SCN}\}$ and $\Sigma^{-}=\{c<\beta_{SCN}\}$,
    \item  $\Gamma^{+}=\{f_{W}>\theta_{W}\}$ and $\Gamma^{-}=\{f_{W}<\theta_{W}\}$.
\end{itemize}
	
These boundaries divide the domain of the model into the following four subregions:
\begin{enumerate}
	\item $\Sigma^{+}\cap\Gamma^{+}=\{c>\beta_{SCN}$ and $f_{W}>\theta_{W}\}$ (wake state with increasing $h$ and high $f_{SCN}$),
	\item $\Sigma^{-}\cap\Gamma^{+}=\{c<\beta_{SCN}$ and $f_{W}>\theta_{W}\}$ (wake state with increasing $h$ and low $f_{SCN}$),
	\item  $\Sigma^{-}\cap\Gamma^{-}=\{c<\beta_{SCN}$ and $f_{W}<\theta_{W}\}$ (sleep state with decreasing $h$ and low $f_{SCN}$),
	\item  $\Sigma^{+}\cap\Gamma^{-}=\{c>\beta_{SCN}$ and $f_{W}<\theta_{W}\}$ (sleep state with decreasing $h$ and high $f_{SCN}$).
\end{enumerate}

In each of these subregions, the model has smooth dynamics dictated by subsets of Equations (\ref{FF_fW}) - (\ref{FF_theta}), while on the boundaries $\Sigma$ and $\Gamma$ dynamics are defined by Filippov's convex method. In Appendix \ref{APX:CHS}, we show that the model flow is transversal across the boundaries of these four subregions and thus, a solution of this piecewise smooth system can be concatenated from trajectories in its four subregions.

\subsection{Bifurcations in the CHS model}
 Using a numerical algorithm by Calvo et al. \cite{Algorithm968} 
 to integrate the CHS model, we numerically computed the bifurcation diagram of the rotation number $\rho$ for $k \in (0,1]$ (Figure \ref{fig:HS_figures}B). The bifurcation diagram maintains a period-adding-like structure, as well as similar trends to those observed for small $\alpha_{SCN}$ values. 
As $k$ was decreased from 1, the $\rho=1$ solution in the CHS model lost stability at $k=0.45$, a similar value as found for $\alpha_{SCN}=0.3$. As suggested by the shrinking of $k$ intervals for solutions with $\rho \in (\frac{1}{2},1)$ for small $\alpha_{SCN}$, the next stable solution detected as $k$ was decreased was $\rho = \frac{1}{2}$. Specifically,  
we did not detect any intermediate $\rho=\frac{q}{p}$ solutions between the $\rho=1$ and $\rho=\frac{1}{2}$ solutions in the CHS model. As $k$ decreased further, there were fewer stable solutions between the $\rho=\frac{1}{2}$ and $\rho=\frac{1}{3}$ solutions in the CHS model than for solutions associated with larger $\alpha_{SCN}$ values. Between the $\rho=\frac{1}{3}$ and $\rho=\frac{1}{4}$ solutions, no intermediate solutions were detected in the CHS model. 

Below, we explain why we detect many or few types of intermediate solutions between certain $\rho$ values using the fast-slow surface associated with the CHS model.
\subsubsection*{Winnowing of entrainment regions 
}

In the CHS model, the smoothly varying $Z$-shaped surface associated with continuous $c(t)$ is split into two connected $Z$-shaped surfaces corresponding to the positive and negative $c(t)$ of the CHS model. By considering model trajectories on this double $Z-$shaped surface, we identify a change in the bifurcation sequence delimiting the $\rho=\frac{1}{2}$ and $\rho=\frac{1}{4}$ stable solutions. We note that the switching boundary $\Sigma$ constrains sleep onset phases for stable $\rho=1$ solutions. Namely, as $k$ is decreased from 1, sleep onset phases remain at $\Phi=0.75$ due to trajectories falling off the upper wake manifold of the $Z-$shaped surface at $\Sigma$ during decreasing circadian drive $c(t)$ (Figure \ref{fig:HS_figures}C). At $k=0.45$ where the $\rho=1$ solution loses stability, the trajectory additionally approaches $\Sigma$ during increasing circadian drive at $h$ values where a transition off the upper manifold is possible (Figure \ref{fig:HS_figures}D). For slightly smaller $k=0.449$, the $\rho=\frac{1}{2}$ solution gains stability in a border collision bifurcation due to sleep onsets occurring at $\Sigma$ for $\Phi=0.25$ and $\Phi=0.75$. 

The absence of solutions with $\rho \in (\frac{1}{2},1)$ is due to the trajectory intersecting $\Sigma$ as $c(t)$ increases and falling off the upper wake manifold at every circadian cycle. A stable solution with $\rho \in (\frac{1}{2},1)$ requires that sleep onset phases slightly shift on successive circadian cycles such that some cycles only have 1 sleep episode and some cycles have 2 sleep episodes. Here, since one sleep onset always occurs at $\Phi=0.75$ and the trajectory during the sleep episode is similar even if the $h$ value at sleep onset is different (Figure \ref{fig:HS_figures}D), trajectories do not shift to avoid falling off the upper wake manifold at $\Phi=0.25$ at the $\Sigma$ boundary.


In contrast, stable solutions with $\rho \in (\frac{1}{3}, \frac{1}{2})$ were obtained.  
The existence of stable $\rho \in (\frac{1}{3}, \frac{1}{2})$ solutions can be understood by considering the sleep onset phases of the multiple sleep episodes in those solutions. When the $\rho=\frac{1}{2}$ solution loses stability at $k=0.28$, sleep onsets occur near $\Phi \approx 0$ near the minimum of $c$ and at $\Phi = 0.75$. Thus, the switching boundary $\Sigma$ constrains the phase of only one of the sleep episodes. For $k$ slightly smaller, i.e. $h$ slightly faster, sleep onsets will occur at earlier phases. For the sleep onset occurring near $\Phi \approx 0$, the phase is not constrained by the boundary $\Sigma$ and can shift such that a third sleep onset may occur in a circadian cycle resulting in a solution with $\rho \in (\frac{1}{3}, \frac{1}{2})$. 
The evolution of sleep onset phases and sleep patterns as the period-adding structure progresses with decreasing $k$ introduces the beginning of the stable $\rho=\frac{1}{3}$ solution with sleep onset phases close to 0, smaller and larger than 0.75.

The $\rho=\frac{1}{3}$ solution loses stability in a border collision bifurcation and directly transitions to the stable $\rho=\frac{1}{4}$ solution similarly to the way in which the $\rho=1$ solution transitions to the $\rho=\frac{1}{2}$ solution.
At the loss of stability of the $\rho=\frac{1}{3}$ solution at $k=0.208$, the three sleep onset phases  occur near the extrema of the circadian drive (i.e. $\Phi \approx 0$ and $\Phi \approx 0.5$) and at $\Phi = 0.75$ at the $\Sigma$ boundary with decreasing $c(t)$. For $k=0.207$, the $\rho=\frac{1}{4}$ solution gains stability as a fourth sleep onset occurs at $\Phi=0.25$ at the $\Sigma$ boundary with increasing $c(t)$. The constraint that the trajectory always intersects $\Sigma$ as $c(t)$ increases, causing the border collision, does not permit the slight shifting of sleep onset phases on successive circadian cycles necessary to result in a solution with $\rho \in (\frac{1}{4}, \frac{1}{3})$. Instead, the trajectory falls off the upper wake manifold at $\Phi=0.25$ on every circadian cycle resulting in the stable $\rho=\frac{1}{4}$ solution. 


\subsubsection*{Understanding the small $\alpha_{SCN} >0$ case}
For small $\alpha_{SCN} > 0$, similar constraints on sleep onset phases near $\Phi = 0.75$ and $0.25$ also explain the shrinking $k$ intervals for stable solutions with $\rho \in (\frac{1}{2}, 1)$ and $\rho \in (\frac{1}{4}, \frac{1}{3})$. Sleep onset phases in these solutions have slightly different values on successive circadian cycles that can result in different numbers of sleep episodes per cycle. For small $\alpha_{SCN}$, the steep $f_{SCN}$ profile similarly constrains sleep onset phases to be near $\Phi = 0.75$ and $\Phi = 0.25$ as observed in the CHS model. This restricts the ability to sustain differences between trajectory orbits on successive circadian cycles and prevents the slight shifts in sleep onset phases necessary for the stability of these solutions. The solutions in these particular $\rho$ intervals are affected because the additional sleep episode occurring in a circadian cycle (the 2nd sleep episode for $\rho \in (\frac{1}{2},1)$ and the 4th sleep episode for $\rho \in (\frac{1}{4},\frac{1}{3})$) occurs at phases near $\Phi = 0.25$. These constraints do not affect the newly obtained sleep onset phase in the stable $\rho=\frac{1}{3}$ solution ($\Phi\approx 0 )$, and therefore more types of solutions with $\rho\in(\frac{1}{3},\frac{1}{2})$ are detected in this regime. The $\rho=\frac{1}{2},\frac{1}{4}$ stable solutions are delimited by a BC-S $\rightarrow$ BC-U $\rightarrow$ SN bifurcation sequence, while the $\rho=\frac{1}{3}$ solution is characterized by the sequence SN $\rightarrow$ BC-U $\rightarrow$ BC-S.

\section{Discussion}

\label{Discussion}

In this study, we analyzed the bifurcations in a high-dimensional, piecewise smooth ODE-based model for a sleep-wake flip-flop model under circadian rhythm modulation. As standard numerical continuation algorithms are not effective for nonsmooth systems, our study highlights how applying multiple techniques that reveal model solution structures and their dependence on parameters can facilitate a full bifurcation analysis. We applied fast-slow decomposition to reveal an underlying $Z-$shaped solution manifold that supported the orbits of stable solutions. Tracking stable orbits on the $Z-$shaped surface as homeostatic sleep drive time constants were varied showed how the profile of the folds of the surface dictated and participated in border collision bifurcations of solutions when solution trajectories made tangent intersections with the folds. For border collision bifurcations of stable solutions, we were able to visualize the tangent intersection of the trajectory with the folds of the $Z-$shaped surface, which informed understanding of border collision bifurcations of unstable solutions. Importantly, knowledge of the $Z-$shaped surface enabled the computation of circle maps for model dynamics, as initial conditions were chosen at the upper fold (saddle-node) curves. The circle maps allowed tracking of fixed point solutions, representing periodic solutions in the model, as parameters varied and identification of saddle-node bifurcations as well as border collision bifurcations of unstable periodic solutions. This holistic approach may be applied to other model systems that defy conventional numerical bifurcation analysis methods.

Our analysis focused on the effects of varying two physiologically-motivated factors that affect timing and duration of sleep episodes: The time constants of the homeostatic sleep drive and the profile of the SCN firing rate. The primary bifurcation sequence delimiting stable solutions as the homeostatic drive time constants were decreased (by decreasing the scaling parameter $k$) was SN $\rightarrow$ BC-U $\rightarrow$ BC-S. This sequence was dictated by the shape of the circle maps which exhibit discontinuities with infinitely increasing slopes on one side and finite slopes on the other side. The SN $\rightarrow$ BC-U sequence reflects the gain and loss of fixed points on the map branch(es) near the infinite slope(s) while the BC-S bifurcation reflects the loss of the stable solution on the other end of the map branch(es). This primary bifurcation sequence was modulated by variation of the profile of the SCN firing rate, through the parameter $\alpha_{SCN}$.

One such modified sequence observed for stable solutions displaying multiple sleep episodes per circadian cycle at smaller values of $k$, for example $\rho \le \frac{1}{3}$ for some $\alpha_{SCN}$ values, was SN $\rightarrow$ BC-U $\rightarrow$ BC-U $\rightarrow$ SN. This sequence occurred due to deformation of circle maps such that the finite slope at the discontinuity was less than -1. Our analysis of the loss of stability of the $\rho=1$ solution as $\alpha_{SCN}$ was decreased provided a clear illustration of how this modification can occur as the finite slope at the discontinuity passes through 1. 

Another modified sequence observed for multiple sleep episode solutions at small $k$ values and larger $\alpha_{SCN}$ values was SN $\rightarrow$ SN. This sequence occurred when circle maps were continuous and tangent intersections of trajectories with the folds of the $Z-$shaped surface did not occur. This was because the homeostatic sleep drive varied sufficiently fast (small $k$) and the SCN firing rate profile varied sufficiently slowly (large $\alpha_{SCN}$).

These bifurcation sequences obtained in the SWFF model are similar to those identified in the classic Two Process model under similar parameter variation. Specifically, Bailey et al. \cite{skeldongaps} performed an analytic bifurcation analysis of the Two Process model as the level of the lower circadian threshold was varied leading to similar transitions between monophasic and polyphasic sleep patterns as obtained when homeostatic time constants are varied. In regimes where the amplitude of the circadian thresholds were sufficiently large, they found SN $\rightarrow$ BC-U $\rightarrow$ BC-S bifurcation sequences delimiting the stable regimes of solutions that followed a period adding sequence. In this regime, the analytically computed circle maps were monotonic and discontinuous with an infinite slope on one side of the gap, similar to the computed circle maps for the SWFF model. In regimes where the circadian threshold amplitudes were small, the circle maps became continuous due to the absence of tangent intersections of the homeostatic sleep process with the circadian thresholds and bifurcation sequences as circadian threshold levels were varied changed to SN $\rightarrow$ SN.  

In our analysis, we found that that the complexity of bifurcations increased for small values of $\alpha_{SCN}$, as may be expected with steeper profiles of the folds of the $Z-$shaped surface. This was especially true for $\rho=\frac{1}{2}$ solutions which displayed intervals of bistability for some $\alpha_{SCN}$ values. Interestingly, the evolution of the two stable $\rho=\frac{1}{2}$ solutions as $k$ decreased changed for different values of $\alpha_{SCN}$. For the highest $\alpha_{SCN}$ values in this region of bistability, the 2nd stable solution gained stability in a saddle node and then a subsequent saddle node bifurcation destroyed the original stable solution, thus the solutions replaced one another as $k$ decreased. For lower $\alpha_{SCN}$ values in the bistability region, the 2nd stable solution gained stability in a saddle node bifurcation and then was destroyed in a border collision (stable) bifurcation, leaving the original solution as the only stable solution. For an intermediate value of $\alpha_{SCN}$, both solutions lost stability at the same $k$ value in a coincident saddle node and border collision (stable) bifurcation, each bifurcation involving one of the stable fixed points. 

Coexistence of stable solutions can occur in piecewise smooth maps with discontinuities \cite{AvrutinBook, Keener1980}. In many such maps showing coexistence of stable solutions, such as bistability, the values of the map branches across the discontinuity cover an overlapping interval. In our maps, there is no overlap of values of the map branches across the discontinuities, instead bistability emerges due to a deformation of the shape of map branches that introduces multiple fixed points. Given the similarity in its dynamics with the SWFF model, the Two Process model may be a good reduced system to analyze this mechanism for bistability. While bistability has not been previously reported in the Two Process model, an analysis in which analogous parameters to those considered here in the SWFF model are varied has not been conducted to our knowledge.

The striking effect of steeper profiles of the SCN firing rate (small $\alpha_{SCN}$) was the winnowing of certain stable solutions, namely the $\rho \in (\frac{1}{2},1)$ and the $\rho \in (\frac{1}{4}, \frac{1}{3})$ solutions. Solutions that persisted in these $\rho$ intervals as $\alpha_{SCN}$ decreased had rotation numbers of the form $n / (n+1)$. A similar winnowing of solutions has been observed in threshold models when the profile of the threshold is a square wave \cite{ArnoldArrhythmias}. By analyzing the CHS model, where the SCN firing rate profile is a square wave, we found that solution winnowing was due to constraints on sleep onset phases near $\Phi = 0.75$ and $\Phi = 0.25$, at the edges of the square wave. Specifically, the steep slope of the SCN profile limited the slight variation in sleep onset phases on successive circadian cycles necessary for $\rho \in (\frac{1}{2},1)$ and the $\rho \in (\frac{1}{4}, \frac{1}{3})$ solutions.

Our work demonstrates that the combined effects of the sleep homeostat and circadian waveform modulate the timing, duration and number of sleep episodes in complex ways. These findings suggest that interindividual differences manifested in the time constants dictating the variation of the homeostatic sleep drive \cite{rusterholz} affect the transition from early childhood sleep schedules that include naps to monophasic nighttime sleep that characterizes adult sleep schedules \cite{JENNI2007321,SALZARULO1992107, Jenni2004}. This transition process could be further modulated by the circadian rhythm. Future work is needed to connect the changes observed in the theoretical context of this simplified model to behavior observed in early childhood development. However, our results suggest a pertinent role of SCN activity profile, which is affected by seasonality and light conditions \cite{VANDERLEEST}, in modulating the effects of homeostatic sleep drive variations.

\appendix
\section*{Appendices}
\section{The structure of the map as $k$ and $\alpha_{SCN}$ vary}
\label{Structure_of_maps} The first return circle maps we have presented are characterized by at least one discontinuity. The discontinuity associated with the bifurcations leading to loss of stability of $(p,q)$ periodic solutions is the one caused by tangencies on the upper saddle-node curve of the fast-slow surface. 
In this appendix we explain that in the regime of rotation number $\rho=\frac{q}{p}$, the $p^{th}$ return map has $p$ discontinuities associated with the appropriate discontinuity of the first return map. Let $\Pi:[0, 1]\rightarrow [0, 1]$ represent the first return map that demonstrates a discontinuity because of a tangency on a saddle-node curve. Then $\Pi([0,1])=[0,1]\backslash I$, where $I$ is some interval. According to our results, the discontinuity occurs close to the peak of the circadian oscillator for $k=1$ and starts shifting towards later phases as $k$ decreases. The interval $I$ that is excluded from the range of the map is associated with the rising phase of the circadian oscillator. As mentioned in \cite{onemap}, during the rising phase of $c$, and hence $f_{SCN}$, the dynamics close to the upper-saddle node curve strongly promote the consolidation of wake. This leads to a horizontal gap in the map, where we force sleep onset by following the eigenvector associated with the unstable manifold at those saddle-node points. This is extensively analyzed in \cite{onemap}. This horizontal gap overlaps in large part with the interval $I$.

The discontinuity in terms of our discrete circle map can be described as follows: Assume that the discontinuity of the map occurs between the points $(x_{1}, \Pi(x_1))$ and $(y_{1}, \Pi(y_1))$. Then $\forall \delta>0$, $\exists \epsilon>0$ such that for $\left|{x-y}\right|<\delta$ for $x, y$ around the discontinuity, then $ \left|\Pi(x)-\Pi(y)\right|>\epsilon$. 
In the regime close to the occurrence of a bifurcation, there is one stable periodic orbit of period $p$ and the map is increasing on either side of the discontinuity. 

If the first return map has a discontinuity between $x_{1}$ and $y_{1}$, with $x_{1}<y_{1}$, then by a "backwards" cobwebbing on the map we can find $x_{2}$, $y_{2}$, so that $x_{1}=\Pi(x_{2})$ and $y_{1}=\Pi(y_{2})$. Since the map is continuous everywhere else and invertible, $x_{2}$ and $y_{2}$ are sufficiently close. We can now repeat the same process, and find $x_{3}$ and $y_{3}$ sufficiently close, so that $x_{2}=\Pi(x_{3})$ and $y_{2}=\Pi(y_{3})$. Finally, when we do this $p-1$ times the first return map contains two sequences $\{x_{2},...,x_{p}\}$ and $\{y_{2},...,y_{p}\}$ that satisfy: $x_{j-1}=\Pi(x_{j})$ and $y_{j-1}=\Pi(y_{j})$ for $j=2,..,p$, respectively. Therefore, the $p^{th}$ iteration of the map, $\Pi^{p}$, has $p-1$ discontinuities, each across $x_{j}, y_{j}$, for $j=2,..,p$. The idea is that the map provides approximate initial conditions $x_{j},y_{j}$ on the same two trajectories that will lead to the discontinuity $x_{1},y_{1}$ after $j$ iterations, for $j=2,..,p$. So, we are approximately looking at the same two trajectories when they crossed the section at a "past" time that will eventually lead them to crossing the section again at phases $x_{1}$ and $y_{1}$ after $j-1$ more times.

Additionally, $\Pi^{p}$ has another discontinuity across $x_{1}$ and $y_{1}$ leading to $p$ total discontinuities. The discontinuity of the first return map persists in the $p^{th}$ iteration, since it takes at least $q$ circadian days for the trajectories across the discontinuity to entrain, i.e to converge to the stable $p^{th}$-order cycle of the map. In other words, if one sleeps at a circadian phase corresponding to the infinite slope branch of the map, it will take a few days to converge to the stable sleep pattern. Hence, $\Pi^{p}$ is divided into $p$ branches bordered by two pairs from elements of the sequences $\{x_{1},...,x_{p}\}$ and $\{y_{1},...,y_{p}\}$. Each branch is increasing and maps a subinterval of $[0, 1]$ to an other interval in the range of $\Pi^{p}$.

Since we compute the maps numerically, it is important to note that the values we obtain from this process might not agree exactly with the computed data points of the map. However, in any case we can predict where the discontinuities in the $p^{th}$ return map will occur within some error and how many should exist coming from the discontinuity of the first return map. 

 Recall that discontinuities due to a tangency at the saddle-node curves are characterized by an infinite slope in the left branch of the map curve and a finite slope in the right branch. Hence, in the appropriate regime of the $k-\alpha_{SCN}$ parameter space and with $k$ decreasing, when we first enter the $\Big\{\rho=\frac{q}{p}\Big\}$-regime, the branches with infinite slope intersect the diagonal $\Phi_{n+p}=\Phi_{n}$ at a saddle-node bifurcation. As $k$ further decreases, the finite slope part of each branch approaches the diagonal which leads to the loss of stability of the $(p,q)$ periodic solution due to a border collision or a saddle-node bifurcation. For higher values of $k$ we observe border collisions, but for lower values the finite parts start curving downwards introducing more saddle-node bifurcations. 

When $k$ is sufficiently small, the map becomes continuous, so only saddle-node bifurcations occur. The homeostatic sleep drive $h$ is fast enough now that it can counteract the wake-promoting effect of $f_{SCN}$. In this transition the vertical gap shrinks and interestingly the length of the horizontal gap also reduces accordingly. For values of $k$ and $\alpha_{SCN}$ that the map is continuous, we see that the bifurcation diagram of the rotation number, $\rho$, becomes more dense and continuous as well.


\section{Computation of the bifurcation diagram of the rotation number $\rho$}
\label{Algorithm_BifDiagramOfRho}

The rotation number $\rho = \frac{q}{p}$ describes the number of circadian days over the number of sleep episodes. To compute the rotation number numerically, we have created an algorithm that detects the repeating pattern of sleep episodes from the model trajectory. 

The algorithm works as follows: For each value of $k$ we simulate the model for 100 days to ensure that it has converged to its stable solution. Simultaneously, we keep track of the sleep onsets and their corresponding preceding circadian minima, i.e., the local minima of the variable $c$, using a detection of the event $f_W=4$ during the decrease of the variable $f_W$. This allows us to compute the circadian phase of each sleep onset in the simulation. 

Starting at the last sleep onset phase recorded, we check the preceding sleep onset phases to detect the previous occurrence of the same phase. Since, all of our results are obtained numerically, we allow for an error of 0.0003 for two phases to be considered "equal". If the length of the subsequence that involves the two "equal" phases and all intermediate sleep onset phases is $p+1$, then the number of sleep episodes in the pattern is defined to be $p$ (this avoids double counting the first/last phase). 

To determine the number of circadian days, $q$, we count the distinct number of circadian minima that correspond to the sleep onset phases of the pattern. 

For some values of $k$, this algorithm did not detect a stable repeating pattern. In that case,  we computed an average $\rho$ as the total number of days divided by the total number of sleep cycles in a simulation lasting 120 days. 

\section{Circadian Hard Switch (CHS) model}
\label{APX:CHS}
In the CHS model, model dynamics are smooth in the four subregions ($\Sigma^{+}\cup\Gamma^{+}$, $\Sigma^{+}\cup\Gamma^{-}$, $\Sigma^{-}\cup\Gamma^{+}$, $\Sigma^{-}\cup\Gamma^{-}$)
and dynamics on the boundaries $\Sigma$ and $\Gamma$ are defined by Filippov's convex method. Specifically, for $\mathbf{X}=\{f_{W},f_{S},f_{SCN},h,c,\theta\}$ we represent the model system as follows:

\[ \dfrac{d\mathbf{X}}{dt}=\begin{cases} 
F_{11}(\mathbf{X}) & \mathbf{X}\in \Sigma^{+}\cap\Gamma^{+}\\
\bar{co}\{F_{11},F_{12}\} & \mathbf{X}\in \Gamma^{+}\cap\Sigma\\
F_{12}(\mathbf{X}) & \mathbf{X}\in \Sigma^{-}\cap\Gamma^{+}\\
\bar{co}\{F_{12},F_{21}\} & \mathbf{X}\in \Sigma^{-}\cap\Gamma\\
F_{21}(\mathbf{X}) & \mathbf{X}\in \Sigma^{-}\cap\Gamma^{-}\\
\bar{co}\{F_{21},F_{22}\} & \mathbf{X}\in \Gamma^{-}\cap\Sigma\\
F_{22}(\mathbf{X}) & \mathbf{X}\in \Sigma^{+}\cap\Gamma^{-}\\
\bar{co}\{F_{22},F_{11}\} & \mathbf{X}\in \Sigma^{+}\cap\Gamma\\
\end{cases}
\]

\noindent where $\bar{co}\{F_{ij},F_{kl}\}=\{F_{ij,kl}=\alpha F_{ij}+(1-\alpha)F_{kl},   \alpha\in[0,1]\}$ is a convex combination of the flows on either side of a switching boundary. The vector fields $F_{ij}(\mathbf{X})$ in the different subregions are defined as follows:
\begin{itemize}
	\item The following differential equations regarding the variables $f_W, f_S, c, \theta$ are true for all $F_{ij}(\mathbf{X})$ in their corresponding subregions:
	
	\begin{equation} \label{HSa_fW}
	\dfrac{df_{W}}{dt}=\dfrac{W_{max}\cdot0.5\cdot\bigg(1+\tanh\Big(\dfrac{g_{scnw}f_{SCN}-g_{sw} f_{S}-\beta_{W}}{\alpha_W}\Big)\bigg)-f_{W}}{\tau_{W}}
	\end{equation}
	\begin{equation} \label{HSa_fS}
	\dfrac{df_{S}}{dt}=\dfrac{S_{max}\cdot0.5\cdot\bigg(1+\tanh\Big(\dfrac{-g_{scns}f_{SCN}-g_{ws} f_{W}-(k_{2}h+k_1)}{\alpha_S}\Big)\bigg)-f_{S}}{\tau_{S}}
	\end{equation}
	\begin{equation}\label{HSa_c}
	\dfrac{dc}{dt}=-\omega\sin{\theta}
	\end{equation}
	\begin{equation}\label{HSa_theta}
	\dfrac{d\theta}{dt}=\omega
	\end{equation}

	\item In $F_{11}(\mathbf{X})$ and $F_{22}(\mathbf{X})$ $f_{SCN}$, where $c>0$, the differential equation of $f_{SCN}$ is: 
	
	\begin{equation} \label{HSa_fSCNs1}
	\dfrac{df_{SCN}}{dt}=\dfrac{SCN_{max}\cdot0.5\cdot\bigg(1+\tanh\Big(\frac{1}{0.7}\Big)\bigg)-f_{SCN}}{\tau_{SCN}}
	\end{equation}
	\begin{equation}\label{HSa_hg1}
	\dfrac{dh}{dt}=\dfrac{h_{max}-h}{\tau_{hw}}
	\end{equation}
	
	On the other hand, in $F_{12}(\mathbf{X})$ and $F_{21}(\mathbf{X})$ the differential equation for $f_{SCN}$ is: 
	
	\begin{equation} \label{HSa_fSCNs2}
	\dfrac{df_{SCN}}{dt}=\dfrac{SCN_{max}\cdot0.5\cdot\bigg(1-\tanh\Big(\frac{1}{0.7}\Big)\bigg)-f_{SCN}}{\tau_{SCN}}
	\end{equation}

	\item Similarly, in $F_{11}(\mathbf{X})$ and $F_{12}(\mathbf{X})$, $h$ is increasing and its differential equation is:
	
	\begin{equation}\label{HSa_hg1}
	\dfrac{dh}{dt}=\dfrac{h_{max}-h}{\tau_{hw}}
	\end{equation}
	
	In $F_{21}(\mathbf{X})$ and $F_{22}(\mathbf{X})$, $h$ is decreasing and its differential equation is:
	
	\begin{equation}\label{HSa_hg2}
	\dfrac{dh}{dt}=\dfrac{h_{min}-h}{\tau_{hs}}
	\end{equation}
	
\end{itemize}

\subsubsection*{Ruling out sliding motions}

To verify that the model flow does not permit the occurrence of sliding motion on the switching boundaries $\Sigma$ or $\Gamma$, we need to determine whether trajectories will always cross $\Sigma$ or $\Gamma$ transversally. 
To this end, let  $g(\mathbf{X})=f_W-\theta_W=0$ define the boundary $\Gamma$ and $v(\mathbf{X})=c-\beta_{SCN}=0$ define the boundary $\Sigma$. Then, $\nabla g=<1,0,0,0,0,0>$ and $\nabla v=<0,0,0,0,1,0>$. For each boundary subregion, the conditions verifying that the flow directions on either side of a switching boundary $\Sigma$ and $\Gamma$ are in the same direction are as follows:

\begin{itemize}
	
	\item On $\Sigma^{+}\cap\Gamma$:
	\begin{align*}
	    \big(\nabla g(\mathbf{X})^{T}\cdot& F_{22}\big)\big(\nabla g(\mathbf{X})^{T}\cdot F_{11}\big)=\\
	    &\Bigg(\dfrac{W_{max}\cdot0.5\cdot\Big(1+\tanh\Big(\dfrac{g_{scnw}f_{SCN}-g_{sw} f_{S}-\beta_{W}}{\alpha_W}\Big)\Big)-f_{W}}{\tau_{W}}\Bigg)^2\geq0.
	\end{align*}

\item	Similarly, on $\Sigma^{-}\cap\Gamma$:
\begin{align*}
  \big(\nabla g(\mathbf{X})^{T}\cdot& F_{12}\big)\big(\nabla g(\mathbf{X})^{T}\cdot F_{21}\big)=\\
  &\Bigg(\dfrac{W_{max}\cdot0.5\cdot\Big(1+\tanh\Big(\dfrac{g_{scnw}f_{SCN}-g_{sw} f_{S}-\beta_{W}}{\alpha_W}\Big)\Big)-f_{W}}{\tau_{W}}\Bigg)^2\geq0  
\end{align*}
	
	\item On $\Sigma\cap\Gamma^{+}$:\\
	 $\big(\nabla v(\mathbf{X})^{T}\cdot F_{1}\big)\big(\nabla v(\mathbf{X})^{T}\cdot F_{2}\big)=(-\omega\sin(\theta))^2\geq0$.

\item	Similarly, on $\Sigma\cap\Gamma^{-}$: \\
	$\big(\nabla v(\mathbf{X})^{T}\cdot F_{3}\big)\big(\nabla v(\mathbf{X})^{T}\cdot F_{4}\big)=(-\omega\sin(\theta))^2\geq0$.

\end{itemize}

All four of these conditions are satisfied indicating that the boundaries cannot be attracting (or repelling) from both sides. This is sufficient to ensure that trajectories transversely cross  each of the switching manifolds $\Sigma$ and $\Gamma$ \cite{DiBernardoBook}. Thus, the possibility of sliding on $\Sigma$ and $\Gamma$ is eliminated, and the representation of the flow on these boundaries as a convex combination of the flow on either side of the boundary is well-defined with an arbitrary choice of $\alpha$. 

\section*{Acknowledgments}
The authors thank Anne Skeldon for helpful discussions and Grace O'Brien for assistance with the 3D surface plots.

\bibliographystyle{plain}

\end{document}